\documentclass[12pt]{amsart}
\usepackage{amscd,amssymb,graphics}
\newtheorem{thm}{Theorem}[subsection]
\newtheorem{corol}[thm]{Corollary}

\theoremstyle{definition}
\newtheorem{defin}[thm]{Definition}

\theoremstyle{remark}
\newtheorem{remark}[thm]{Remark}
\newtheorem{example}[thm]{Example}
\newtheorem{examples}[thm]{Examples}
\newtheorem{question}[thm]{Question}
\newtheorem{problem}[thm]{Problem}
\newtheorem{conj}[thm]{Conjecture}
\numberwithin{equation}{section}

\newcommand{\abs}[1]{\lvert#1\rvert}
\def\norm#1{\left\Vert#1\right\Vert}
\def\U {{\mathbb U}}
\def\R {{\mathbb R}}
\def\I {{\mathbb I}}
\def\K {{\mathbb K}}
\def\C {{\mathbb C}}
\def\N{{\mathbb N}}
\def\e{{\epsilon}}
\def\Z {{\mathbb Z}}
\def\s{{\mathbb S}}
\def\H{{\mathcal H}}
\def\Homeo{{\mathrm{Homeo}\,}}
\def\Spec{{\mathrm{Spec}\,}}
\def\Iso{{\mathrm{Iso}\,}}
\def\Aut{{\mathrm{Aut}\,}}

\begin{document}

\title[Topological groups]{Topological groups: 
where to from here?}

\author[V. Pestov]{Vladimir Pestov}
\address{School of Mathematical and Computing Sciences,
Victoria University of Wellington, P.O. Box 600, Wellington,
New Zealand}
\email{vova@mcs.vuw.ac.nz}
\urladdr{http://www.mcs.vuw.ac.nz/$^\sim$vova}

%    General info

\date{September 30, 2000}


\begin{abstract} This is an account of 
one man's view of the current perspective
of theory of topological groups. 
We survey some recent developments which are, from our viewpoint,
indicative of the future directions, 
concentrating on actions of topological groups
on compacta, embeddings of topological groups, free
topological groups, and `massive' groups (such as groups of homeomorphisms
of compacta and groups of isometries of various metric spaces).
\end{abstract}

\maketitle

\setcounter{tocdepth}{2}
%{\small
\tableofcontents
%}

%\setcounter{tocdepth}{2}
%\tableofcontents

\section{Introduction}
\subsection{Motivation}
For a randomly selected mathematician outside of the field of
general
topology --- or, to be more precise, `general topological algebra' ---
the words `topological group' most probably sound synonymous with
`locally compact group.' Indeed, the depth, beauty and importance of
theory of locally compact groups, in particular
representation theory, abstract harmonic analysis, duality theory
and structure theory, are overwhelming, while the 
richness of links with other areas of mathematics, physics, chemistry,
computing and other sciences is hard to match. 
A natural question to ask is therefore:
\\

$\bullet$ {\it Is there life beyond local compactness?}
\\

In the present audience the question sounds hollow rhetoric as
the answer is known to be in the affirmative to just about everyone.
Therefore, I will sharpen it up:
\\

$\bullet$ {\it Does there exist a potential for a 
theory of {\rm (}some classes of{\rm )} 
non-locally compact groups of
a comparable depth to that of theory of locally compact groups?}
\\

The answer to this question is much less obvious, as probably only
very few
people would rush into betting on a definite `yes' or `no.'
There is a not uncommonly held opinion 
(which some of us would dispute) that
the most recent result in the theory of topological groups that really mattered
for the rest of mathematics was the solution to Hilbert's 
Fifth Problem 
\cite{MZ}; can we presently discern where the next result about 
topological groups to have an equally resounding impact could
emerge from? 
As a matter of fact, every professional working in the field is guided 
by a vision of his/her own, and the totality
of views coming from different researchers would span
a remarkably diverse range of opinions. 
This paper offers such an opinion belonging to 
the present author. Accordingly, the 
article does not attempt to be in any way 
comprehensive and compete with
such substantial surveys on topological groups 
as \cite{CHR,Tk4}, or else \cite{Arh1,Arh2}. 

One obvious approach to exploring the general 
question stated above is to try and directly extend concepts
and results from the locally compact case to more general classes of groups.
In this way, one surrounds the class of locally compact groups with
a larger `halo' formed by topological groups that inherit some or other
features of locally compact groups.
One of such species of topological groups thriving in
the penumbral shadows of local compactness is the class of nuclear (abelian)
groups advertised by Banaszczyk \cite{Ba2}, 
whose {\it raison d'\^etre} is essentially 
testing the limits of Pontryagin--van Kampen duality, as well as
of the entire
body of closely related structural results such as Glicksberg's 
theorem \cite{BMP}. Another such
species is the class of pseudocompact groups which is very popular nowadays 
largely due to the consistent 
efforts of Wis Comfort over the past few decades \cite{CR, Comf}.
There can hardly be
a better recipe for achieving a deeper insight into the nature of
(locally) compact groups than walking a few steps
away and having a good look from the outside! 
As a matter of fact, this approach can lead to
fruitful insights into the structure and properties of topological
groups of very general nature. One would expect essential further
progress achieved in this direction in the future.
Nevertheless, one meets certain limitations on this way, since
most of tools making the theory of LC groups a
success are intrinsic to local compactness -- like the
left-invariant locally finite finitely additive Borel measure, 
positive on non-void open sets \cite{Al}.

This makes an alternative approach unavoidable: to isolate new 
classes of topological groups of importance on their own right
with a view towards understanding their structure and
properties. 
Such large (`esentially non-locally compact') 
topological groups originate in many different contexts:
in set-theoretic and smooth topology, ergodic theory, representation
theory, functional analysis, and topological dynamics, to name a few. 
Below we will exhibit examples of
such groups and show that some of their 
properties have no analogue in the locally compact case.

\subsection{What is included in the article}
The paper is loosely `coordinatized' by the following four notions.
\par
{\it 1. Embeddings.}
\par
The following question has stimulated investigations in the
theory of topological groups over many years.
Let $\mathcal P$ be a non-empty class of topological groups
(possibly consisting just of a single group). When
is a given topological group $G$ isomorphic to a subgroup of
a group from $\mathcal P$?
\par
2. {\it Actions.}
\par
The situation where a topological group $G$ acts continuously on a
topological space $X$ emerges very often in disparate 
contexts throughout mathematics. 
If $X$ is compact, then the triple formed
by $X$, $G$, and the action of $G$ on $X$ forms, formally speaking,
the object of study of abstract topological dynamics. However,
the accents in topological dynamics are put on different concepts ---
having more to do with the structure and properties of orbits --- from
those we are interested in. Our emphasis will be also different from
theory of $G$-spaces, where the main attention is paid to the
phase space $X$. Therefore, we prefer to talk simply of
actions.

3. {\it `Massive,' or `large,' groups.}

It is impossible to give a 
formal definition of this concept, yet `large' groups are easy to
recognize. Examples include
the group of all self-homeomorphisms of a 
(sufficiently homogeneous) compact
or locally compact space equipped with the compact-open topology
(or related topologies),
the full unitary group of an infinite-dimensional Hilbert space,
typically with the strong operator topology, the group of
measure-preserving transformations of a Lebesgue measure space with
the weak topology, and so forth. The role played by such groups
in mathematics is fundamental, yet there is no coherent unifying
theory in existence for the time being. 

4. {\it Free topological groups.}

Free topological groups, introduced by Markov in 1941 along with their
closest counterparts such as free abelian topological groups and free
locally convex spaces, served as an inspiration for the concept
of a universal arrow to a functor introduced by Pierre Samuel.
Free topological groups remain a very useful source of examples
and building blocks in general topological group theory.
Apart from that, these objects have never enjoyed much 
popularity and are often perceived as exotic. However, such
objects, properly disguised, often resurface in other areas of
mathematics, meaning that the accumulated
expertise of topological group theorists is very probably applicable to
problems that at the first sight have nothing to do with topological
algebra.

\smallskip
Each of the above four concepts will feature below in three different
roles: as a tool, as an object of study, and as a link between
`general topological algebra' and other disciplines of mathematics.
Notice that all of the above are closely intertwined and their role
in our presentation is that of a
coordinate system rather than a linear index.

\subsection{What is left out} It is hardly surprising that the
vast majority of research directions pertinent to `large' topological
groups is left untouched by the present paper.
As the first such omission, we want to mention descriptive set theory,
or more exactly descriptive theory of group actions,
where non-locally compact Polish (that is, completely metrizable)
topological groups appear very naturally. A topological algebraist must
certainly keep one eye on further developments in the area, which
provides both a motivation and a guidance 
for the future general theory of non-locally compact 
topological groups. The present
author does not feel qualified to touch upon this subject, 
but fortunately there are excellent references available, e.g.
\cite{BK,Kech}.

Another omission, but this time quite purposeful, 
is theory of infinite-dimensional Lie groups.
Indeed, many concrete large topological groups support a natural
structure of ($C^\infty$ or sometimes even analytic)
Lie groups modelled over locally convex spaces of infinite dimension.
While such a theory for Banach--Lie groups goes back to the 30's
and is well established \cite{bou,HM}, 
there are several competing versions of
infinite-dimensional Lie theory beyond the Banach--Lie case
\cite{KYMO,Mil,Les} (see also \cite{Neeb}, which was referee's
suggestion).
Already for Lie groups modelled over nuclear Fr\'echet spaces the
theory meets some difficulties of fundamental nature, the most
upsetting of which is the apparent need to explicitely
require the existence of exponential map: it is still an open problem
whether or not a $C^\infty$ Fr\'echet manifold equipped with smooth
group operations admits an exponential map from the tangent 
Lie algebra at the identity ({\it loco citato}).

However, the reason why infinite-dimensional groups are not featured
in this essay is of a different kind: for all we know, they cannot be
considered as mere {\it topological groups,} but rather as
{\it topological groups with additional structure} --- unlike in
the finite dimensional case, where the classical Montgomery--Zippin
theorem allows one to identify Lie groups with 
locally Euclidean topological groups! No analogue of
the Montromery--Zippin theorem 
is known even for Banach--Lie groups. Moreover, such an analogue
is impossible to state in terms of
topology alone, in view of the following result by Keesling \cite{Kee}:
every separable metrizable topological group embeds as a topological
subgroup into a topological group homeomorphic to the separable
Hilbert space $l_2$. (Now recall that for every Banach--Lie 
group the restrictions of the two natural uniformities to a
suitable neighbourhood of identity coincide, and this property,
while inherited by topological subgroups, is typically 
absent in Polish topological groups.)

Finally, it is beyond reasonable doubt 
that topological {\it semigroups} will
play an ever increasing role in theory of `large' topological groups.
Semigroups do appear in this survey on a few occasions, 
but present author does not feel
sufficiently competent in the area to offer topological semigroups
the prominent place they richly deserve. A couple of references 
on the subject are \cite{BJM,CHR}.

\subsubsection*{\bf Acknowledgements}
The author expresses his gratitude to the Organizers of the
14-th Summer Conference on General Topology and its Applications for
selecting him as an Invited Speaker and thus offering an 
opportunity to put in order his thoughts about topological groups.
This article is a (somewhat extended) write-up of the actual
lecture.

Thanks also go to Vladimir Uspenskij, Michael Megrelishvili, and
the anonymous referee for their constructive criticisms
of the initial version of this paper (dated 29.10.99), and to Alekos Kechris
for useful remarks.

Some of the research the article is based upon was supported by the
Victoria University of Wellington Research Development Fund.

\section{Actions and representations}

\subsection{Some basics}
We will be only considering Hausdorff topological 
groups, Tychonoff
topological spaces, and separated uniform spaces. 
A topological group $G$ {\it acts} on a 
topological space $X$ if there is a continuous
mapping $\tau\colon G\times X\to X$, called an {\it action,}
where the image of a pair $(g,x)\in G\times X$ is usually denoted either by
$\tau_gx$ or simply by $g\cdot x$, having the properties that
$g\cdot (h\cdot x)=(gh)\cdot x$ and
$e_G\cdot x=x$ for every $g,h\in G$ and $x\in X$. Here and
in the sequel, $e=e_G$ denotes the
identity element of the acting group. 
The entire triple ${\mathfrak X}=(X,G,\tau)$ is called
an (abstract) {\it topological dynamical system,} while $X$ together
with the action $\tau$ is referred to as a $G${\it -space.}
The triple $\mathfrak X$ is also known under the name a
{\it topological transformation group.}

If the space
$X$ is compact, a continuous action of a topological group $G$
on $X$ can be identified with a topological group
homomorphism $G\to\Homeo_c X$, where the subscript $`c'$ denotes the
compact-open topology on the homeomorphism group. Namely, to every
element $g\in G$ one associates the mapping 
\[X\ni x\mapsto g\cdot x\in X,\]
which is a self-homeomorphism of $X$ (with the inverse
$x\mapsto g^{-1}\cdot x$). Of course, this mapping makes sense
for a non-compact space $X$ as well, and is called 
a {\it motion}, or a {\it translation}. However, only in the case of
compact $X$ one can always prove that the homomorphism
$G\to\Homeo_c X$ is continuous. 
Conversely, every continuous homomorphism from a
topological group $G$ to the full homeomorphism group 
$\Homeo_cX$ of a compact space $X$
determines in a unique way a continuous action of $G$ on $X$.
For general topological spaces $X$ the
group of self-homeomorphisms $\Homeo X$ no longer 
supports an apparent group
topology making such a convenient identification possible, though
there are important exceptions.

Let us single out two particularly significant classes of actions.
An action is called {\it effective}
if every element $g\in G$, different from the identity, acts in a
non-trivial fashion on {\it some} element of the phase space $X$,
that is,  $\forall g\neq e$, $\exists x$, $g\cdot x\neq x$.
An action is called {\it free} if every element $g$ different from
the identity acts in a non-trivial way on {\it every} element of $X$,
that is, $\forall g\neq e$, $\forall x$, $g\cdot x\neq x$.

\begin{examples} 1. Every topological group acts on itself freely
by left translations through $g\cdot x=gx$, $g,x\in G$. \\
2. The action of a topopological group $G$ on itself by conjugations,
\[g\cdot x= g^{-1}xg, ~~ g,x\in G,\]
is effective if and only if the centre of $G$ is trivial.\\
3. The canonical action of
the group $\Homeo_c\s^1$ on the circle $\s^1$ is effective but not free.
\end{examples}

The topological space $X$ (the {\it phase space}) 
can also support additional structures
of various sort. If $X=(X,\rho)$ is a metric space, then an action of a
topological group $G$ on $X$ is said to be {\it isometric}, or an
action {\it by isometries,} if every motion $x\mapsto g\cdot x$,
$g\in G$, is an isometry of $X$ onto itself. In this case, the
continuity of the action is equivalent to the continuity of the
associated homomorphism from $G$ to the group of isometries
$\Iso(X)$ equipped with the pointwise topology (that is, one
inherited from $X^X$, or else from $C_p(X,X)$). 
Moreover, the group $\Iso(X)$ with the pointwise topology forms
a Hausdorff topological group. Both statements are very easy to
verify directly.

If $X$ is a Banach space --
more precisely, a complete normed linear space -- and $G$ acts on $X$
continuously, then the action $\tau$ is called a {\it continuous
representation of $G$ in $E$} if every motion $E\ni x\mapsto g\cdot x\in E$,
$g\in G$, is a linear operator.
Effective actions are known in this context
as {\it faithful representations}. Every such action has
zero element as the fixed point and therefore cannot be free.
One says that the action $\tau$ on a Banach space $E$
is {\it by isometries} if every
homeomorphism of the form $x\mapsto g\cdot x$, $g\in G$, is a linear
isometric transformation of $X$. In such a case, the action of $G$
on $X$ is often referred to as a representation of $G$ in $X$ {\it by
isometries.} For such actions instead of continuous 
representations one normally speaks of
{\it strongly continuous} representations. To explain this peculiarity
in terminology, recall that the
{\it strong operator topology} on the space ${\mathcal L}(E,E)$
of all bounded linear operators on a locally convex space $E$
is the topology of pointwise convergence, that is,
one induced from the standard embedding
${\mathcal L}(E,E)\hookrightarrow E^E$, where the latter space carries
the Tychonoff product topology. 
The restriction of the strong operator topology to the group
$\Iso(E)_s$ therefore coincides with the topology of pointwise
convergence and,
as we noted above, an isometric representation of $G$ on $E$
is continuous if and only if the associated homomorphism
$G\to\Iso(X)$ is continuous with respect to the strong operator topology
on the latter group.

Finally, if $X=\H$ is a (complex)
Hilbert space, then a representation $\tau$ of a
group $G$ in $\H$ is called {\it unitary} if every motion
$\tau_g\colon\H\ni x\mapsto g\cdot x\in\H$ is a unitary operator:
$(\tau_gx,y)=(x,\tau_{g^{-1}}y)$ for all $x,y\in\H$. It is a particular
case of a representation by isometries, and in fact namely strongly continuous
unitary representations of topological groups are of overwhelming
importance. The collection of all unitary operators on a Hilbert
space $\H$ forms a group, called the {\it full unitary group} of $\H$
and denoted by $U(\H)$. The subscript `s' will denote the strong
operator topology, and the topological group $U(\H)_s$
(where $\dim\H=\infty$) is one of the
most important `massive' groups. If $\H$ has finite
dimension $n$, then $U(\H)=U(n)$ is the group of $n\times n$ unitary
matrices with the natural compact topology.

\subsection{Teleman's theorem}
The following fundamental result asserts, roughly speaking,
that effective actions of
topological groups are sufficiently common.

\begin{thm}[Silviu Teleman, 1957, \cite{Te}] Every Hausdorff 
topological group $G$ acts effectively \\
($\ast$) on a Banach space by isometries; \\
($\ast\ast$) on a compact space.
\label{teleman}
\end{thm}

Before proceeding to the proof, we will say a few words about
uniformities on topological groups.
The difference between
the two standard uniform structures on a topological group
turns out to be surprisingly important and, in particular,
provides the clue to some of the
subtlest results in theory of topological groups and their actions
known to date,
which is why taking special care of terminology and notation 
is important. If $G$ is a
topological group, then the {\it right equicontinuous}
(in short, {\it right e.c.}) uniform structure
on $G$, which we will denote by ${\mathcal U}_\Rsh(G)$, has basic
entourages of the diagonal of the form
\[V_\Rsh:=\{(g,h)\in G\times G\colon gh^{-1}\in V\},\]
as $V$ runs over a neighbourhood basis of $G$. The uniformity 
${\mathcal U}_\Rsh(G)$ is right equicontinuous in the sense that 
the family of all right translations $R_g$, $g\in G$, where
$R_g(x):=xg$, is uniformly equicontinuous as a family of mappings from
the uniform space $(G,{\mathcal U}_\Rsh(G))$ to itself. [The 
basic entourages $V_\Rsh$ are invariant with respect to the 
action of $G$ on $G\times G$ on the right:
\[(h,f)\cdot g:=(hg,fg).]\] 
The right equicontinuous uniformity is called 
{\it right uniformity} by some
authors and {\it left uniformity} by others (consult our paper
\cite{P19} for references to the both kinds of usage). Since the
word {\it right equicontinuous} is not only unambiguous but fully
descriptive, we will adopt it in this article. In a similar way, one
can consider the {\it left equicontinuous} uniformity (suggested notation:
${\mathcal U}_\Lsh(G)$), determined by the basic neighbourhoods of the
form 
\[V_\Lsh:=\{(g,h)\in G\times G\colon g^{-1}h\in V\}.\]

Denote by $C^b_\Rsh(G)$ the vector space formed by
all bounded (real or complex-valued)
functions on $G$ that are uniformly continuous with respect to the
right e.c. uniform structure:
\[\forall \e>0, \exists V\ni e, xy^{-1}\in V \Rightarrow
\vert f(x)-f(y)\vert <\e.\]
If equipped with the supremum norm,
\[\norm f :=\sup_{x\in G}\vert f(x)\vert,\] 
$C^b_\Rsh(G)$ becomes a Banach space.

In a similar way, one can define the Banach space $C^b_\Lsh(G)$
of all bounded functions that are uniformly continuous with
respect to the left e.c. structure on $G$.

\begin{remark}
If the two uniformities on a group $G$ coincide, then one
says that $G$ is a {\it SIN group} --- from 
{\it S}mall {\it I}nvariant 
{\it N}eighbourhoods ---
or else a {\it balanced group}. The SIN property of a topological
group $G$ is equivalent to the existence of a
neighbourhood basis at $e_G$ formed by invariant sets $V$, that is,
such that
$gVg^{-1}=V$ for all $g\in G$. For example, every compact and every
abelian topological group are SIN. Another 
prominent example of a SIN group
is the full unitary group $U(\H)$ of a Hilbert space $\H$ equipped
with the {\it uniform} (not strong!) operator topology, that is,
the topology induced by
the natural embedding $U(\H)\subset C^b(B,\H)$, where $B$ is
the closed unit ball in $\H$ and the space of bounded
continuous functions is equipped with the supremum norm.

Clearly, for every SIN group $G$
the two function spaces are 
identical: 
\[C^b_\Rsh(G)=C^b_\Lsh(G).\] 
Rather astonishingly,
it remains unknown to date if the converse is true!
A similar question can be asked about spaces of functions that
are not necessarily bounded: suppose $C_\Rsh(G)=C_\Lsh(G)$,
is then $G$ necessarily a SIN group? This question was first asked by 
Itzkowitz \cite{Itz}, and
a positive answer has been since obtained for a number of particular
(and quite disparate) classes of topological groups.
(See. for instance, \cite{MNP} and references contained in the paper.)
It remains also unknown if the two
versions of Itzkowitz's question (bounded and unbounded ones) 
are equivalent between themselves. A comprehensive survey of what is
known of the problem to date, written by Itzkowitz himself, appears
in this Volume of Topology Proceedings.
\end{remark}

Now --- back to the proof of Teleman's result, which is 
remarkably straightforward.

($\ast$) 
The group $G$ acts on the Banach space
$C^b_\Rsh(G)$ by left translations:
\[(g,f)\mapsto L_g(f), ~~ L_g(f)(x):= f(gx).\]
The action $L\colon G\times E\to E$ is easily verified to be a continuous
and effective action and each operator $L_g\colon E\to E$ 
is a linear isometry. The first claim is established.

\begin{remark} Observe that in the above proof the group $G$
acts by {\it left} translations on the space of functions that are
uniformly continuous with respect to the {\it right} equicontinuous
uniform structure! As a particularly instructive exercise, the 
reader is advised to try and find where the proof breaks down if one
considers the (perfectly well-defined) action of $G$
by left translations on the space $C^b_\Lsh(G)$.
\end{remark}

($\ast\ast$) 
The group $G$ acts in a natural way on the dual space of
$E=C^b_\Rsh(G)$ as follows:
$(g\cdot\phi)(f):=\phi(g^{-1}\cdot f)$, where $\phi\colon E\to\K$
is a linear functional (an element of $E'$), $g\in G$, and
$f\in E=C^b_\Rsh(G)$.
This action is a representation of $G$ in the dual Banach space
by isometries (the so-called {\it contragredient representation}
to the left regular representation $L$) and consequently
can be restricted to the unit ball, $B$, of $E'$.
It is now an easy exercise to verify that the action of $G$ on $B$
is continuous if $B$ is equipped with the 
(always compact) weak$^\ast$ topology (that is, the weakest topology
making each evaluation
mapping of the form $E'\ni \phi\mapsto \phi(f)$, $f\in E$ continuous).

In fact, one can verify somewhat
more: that the continuous monomorphism
$G\hookrightarrow \Homeo_c (B_w)$ is an embedding of topological 
subgroups. Here $B_w$ is the dual ball with the weak$^\ast$ topology.

\begin{corol}
Every topological group $G$ is topologically isomorphic with a
subgroup of the group of homeomorphisms of a suitable compact space
equipped with the compact-open topology. \qed
\end{corol}

One can reformulate the above statement by saying that every
topological group acts {\it topologically effectively} on an
appropriate compact space.

\begin{remark}
One of the most interesting and best known 
open problems related to actions of
topological groups on compacta is the
{\it Hilbert--Smith conjecture,} which can be given the following
equivalent reformulation: if a zero-dimensional
compact group $G$ acts effectively
on a finite-dimensional topological manifold $X$, then $G$ is
necessarily finite. See the related papers \cite{New, Y, RS, Mart}.
\end{remark}

\subsubsection*{A universal group with countable base}
Now we will proceed to a particularly nice
application of Teleman's construction.
In mid-30's Ulam asked in the Scottish Book (cf. Problem 103 in
\cite{Mau}) if there existed a {\it universal separable topological group,}
that is, a separable group $G$ containing an isomorphic copy of every
other separable topological group. As it was discovered
rather soon, the negative answer follows from simple set-theoretic
considerations: indeed, there are more pairwise non-isomorphic
separable groups than there are different subgroups in any single separable
group ({\it ibid.}) In his comments to Ulam's problem,
Kallman had suggested several `corrected' versions
of the same question, one of them being as follows: does there exist a
universal {\it second-countable} topological group? Apparently quite
independently, the same question was promoted by Arhangel'ski\u\i\
\cite{Arh1,Arh2}.
The answer is
contained in the following result, which exemplifies an interplay between
embeddings and actions (as well as massive groups such as
$\Homeo(\I^\omega)$ undoubtedly is).

\begin{thm}
[Uspenskij, 1985, \cite{U2}] $\Homeo_c(\I^\omega)$ 
is a universal second-countable topological group.
\end{thm}

\begin{proof}
Let $G$ be a second-countable group. Then the Banach space 
$C^b_\Rsh(G)$ conains a $G$-invariant separable Banach subspace $E$
whose elements separate points and closed subsets of $G$.
The unit ball $B_w$ of the dual space to $E$
in the weak$^\ast$ topology, being a convex compact 
subset of the separable Fr\'echet
(= completely metrizable locally convex)
space $E'_{w^\ast}$, is homeomorphic to the 
Hilbert cube $\I^\omega$ by force of Keller's theorem (cf. \cite{BP}). 
Consequently,
\[G\hookrightarrow\Homeo_c(B_w)\cong\Homeo_c(\I^\omega).\]
\end{proof}

\begin{remark}
A more accurate rendering of the same idea 
shows that the pair $(\Homeo(\I^\omega),\I^\omega)$ 
forms a {\it universal second-countable topological transformation
group:} every compact $G$-space $X$, where both $G$ and $X$ are 
second-countable, embeds into $(\Homeo(\I^\omega),\I^\omega)$ in a clear sense.
This result seems to have been 
obtained by Megrelishvili independently from Uspenskij, 
only to be published a decade later in \cite{Meg5}.
\end{remark}

\subsection{Urysohn metric spaces and their groups of isometries}
We are going to steer towards another example of
a universal group with countable base. 
A metric space $M$ is a (generalized) {\it Urysohn space} if
for every finite metric space $X$ and every finite
subspace $Y$, every isometric embedding $Y\hookrightarrow M$ extends to
an isometric embedding of $X$ into $M$. 
(Cf. \cite{Ur}, \cite{Kat}, or \cite{Gro}, 3.11$_+$.) 
Every separable Urysohn metric space $M$ 
contains an isometric copy of every other separable metric
space, and if $M$ is in addition complete, it is unique up to
isometry; we will denote it by $\U$.
(A proof of the latter statement consists of shuttling between the
two spaces and building up a recursive sequence
of extensions using increasing chains of finite subspaces
with everywhere dense union chosen in each space.)
The groups of the form $\mathrm{Iso}\,(M)_s$
provide yet another series of examples of `massive' topological groups.

A metric space $X$ is called $n${\it -homogeneous}, where $n$ is 
a natural number, if every isometry between two subspaces of $X$
containing at most $n$ elements each extends to an isometry of
$X$ onto itself. If $X$ is $n$-homogeneous for every natural $n$,
then it is said to be $\omega${\it -homogeneous}.
The complete separable Urysohn space $\U$ is $\omega$-homogeneous
and moreover enjoys the stronger property: 
every isometry between two compact subspaces of $X$ extends to
an isometry of $X$ onto itself. At the same time, 
non-separable Urysohn metric spaces need not have this property.

The following construction of a Urysohn metric space extension of
a given metric space was suggested by Kat\u etov \cite{Kat}.
Let $X$ be a metric space and let $Y\subseteq X$ be a metric subspace. 
Let us say, following \cite{Kat,U3}, 
that a 1-Lipschitz real-valued function $f$ on a 
$X$ is {\it supported on,}
or else {\it controlled by,} $Y$, 
if for every $x\in X$
\[f(x)= \inf\{\rho(x,y) + f(y) \colon y\in Y\}.\]
In other words, $f$ is the largest among all
1-Lipschitz functions on $X$
having the prescribed restriction to $Y$. As an example, every distance
function $x\mapsto\rho(x,x_0)$ 
from a point $x_0$ is controlled by a singleton, $\{x_0\}$.
Remark that every 1-Lipschitz function on $Y$ extends in a unique 
way to a 1-Lipschitz function on all of $X$ controlled by $Y$.

Denote by $X^\dag$ the
collection of all 1-Lipschitz functions $f\colon X\to\R$ 
that are controlled by finite subspaces of $X$ (depending on the
function). If equipped with the supremum metric,
\[d_{X^\dag}(f,g):=\sup_{x\in X}\vert f(x)-f(y)\vert,\]
$X^\dag$ becomes a metric space of the same weight as $X$. 
Moreover, $X$ isometrically embeds into $X^\dag$
in a natural way:
\[X\ni x\mapsto [d_x\colon X\ni y\mapsto d(x,y)\in\R]\in X^\dag.\]
The embedding $X\hookrightarrow X^\dag$ has a much stronger
property than being just isometric: 
\smallskip

(A) whenever $Y$ is a finite metric subspace of $X$
and $Y'=Y\cup\{x^\ast\}$ is an arbitrary 
one-point metric extension of
$Y$, the metric space $Y'$ is isometric to the space
$Y\cup\{f\}$ for some $f\in X^\dag$ under the identification
$x^\ast\leftrightarrow f$. 
\smallskip

Indeed, the distance function $f\colon x\mapsto d_X(x,x^\ast)$ is in
$X^\dag$, and the metric spaces $Y'$ and $Y\cup\{f\}$ are
isometric.

Apply the Kat\u etov extension $X\mapsto X^\dag$ to an arbitrary
metric space $X$ in a recursive fashion $\omega$ times, and
denote by $\widetilde X$ the completion of the metric space union
of all finite iterations:
\[\cup_{n\in\N}X^{\dag\dag\dag\cdots\dag ~ \mbox{\tiny ($n$ times)}}.\]
The metric space $\widetilde X$ is generalized Urysohn. In particular,
if $X$ is separable, then $\widetilde X$
is isometric to a complete separable Urysohn space $\U$. 

The following theorem also appears 
as an exercise in Gromov's book 
\cite{Gro}, Ch. 3$\frac 12$, published in 1999.

\begin{thm}
[Uspenskij, 1990, \cite{U1}] 
The topological group $\mathrm{Iso}\,(\U)_s$
is a universal second-countable group.
\label{ury}
\end{thm}

\begin{proof}
The proof is based on the following remarkable
property of the functor $X\mapsto X^\dag$, first noticed
and put to use by Uspenskij:
every action of a group $G$ on the space $X$ by isometries 
extends in a canonical way to an action of $G$ on $X^\dag$ (by
left translations), and if the original action of $G$ on $X$ was
continuous, so will be the extended action of $G$ on $X^\dag$.
(To better appreciate the usefulness of Kat\u etov
functions, notice that in general the same action of $G$ 
on the space of {\it all} 1-Lipschitz functions on $X$
need {\it not} be continuous!) 
Now it is rather evident that
every continuous action of a topological group $G$ on
a metric space $X$ by isometries
extends to a continuous action of $G$ on $\widetilde X$ be isometries
in a canonical sort of way. If the original action on $X$ determined an
embedding of topological groups $G\hookrightarrow\Iso(X)_s$,
then clearly $G$ is a topological subgroup of
$\Iso(\widetilde X)_s$ as well.

If $G$ is a separable topological group, then one can start with
any separable metric space $X$ whose group of isometries $\Iso(X)_s$
contains $G$ as a topological subgroup to obtain the desired conclusion.
(For example, $X=E$ as in the proof of Teleman's theorem, or simply
$X=G$ itself equipped with a right-invariant metric generating the
topology.) The space $\widetilde X$ is then separable and
isometric to $\U$ as in the statement of the theorem.
\end{proof}

A more complicated argument, due to the same author, 
establishes the following result.

\begin{thm} 
[Uspenskij, 1998, \cite{U3}] Every topological group
$G$ embeds as a topological subgroup into the isometry group of
a suitable generalized Urysohn space $M$ that is $\omega$-homogeneous
and has the same weight as $G$.
\label{uspee}
\end{thm}

The proof resembles that of Theorem \ref{ury}. However,
in order to achieve $\omega$-homogeneity of the
union space, one has to alternate between
the Kat\u etov metric extension $X^\dag$
and the `equivariant homogenization' extension, $H(\cdot)$. This
extension, forming the nontrivial technical core of the proof, is
described in the following result.

\begin{thm}[Uspenskij \cite{U3}]
Every metric space $X$ embeds, as a metric subspace, 
into an $\omega$-homo\-geneous metric
space $H(X)$ of the same weight as $X$ in such a way that there is a
continuous group homomorphism $e\colon \Iso(X)\to\Iso(H(X))$ with
the property that for each $g\in\Iso(X)$, $e(g)\vert_X=g$. \qed
\label{homo}
\end{thm}

However, since there is no apparent reason for
non-separable complete Urysohn spaces to be unique up to
isometry, the above result cannot be used in order to
answer the following open question:

\begin{question} Let $\tau>\aleph_0$ be a cardinal.
Does there exist a universal topological group of 
weight $\tau$?
\end{question}

\begin{remark} A recent result by Shkarin \cite{Shk} states
that there exists a universal {\it abelian} second-countable
topological group. 
\end{remark}

Returning to the groups of isometries of generalized Urysohn spaces,
we want to mention the following result, whose proof is based on
the lower compactification of the groups $\Iso(U)$.
(Cf. subsection \ref{compacti}.) Recall that a topological group $G$
is called {\it minimal} if it admits no strictly coarser Hausdorff
group topology, and {\it topologically simple} if $G$ contains no
closed normal subgroups other than $G$ itself and $\{e_G\}$.

\begin{thm}[Uspenskij, \cite{U3}]
Let $U$ be an $\omega$-homogeneous generalized Urysohn metric space.
Then the group $\Iso(U)$ is minimal and topologically simple. \qed
\end{thm}

\begin{corol}[{\it loco citato}]
Every topological group $G$ embeds, as a topological subgroup, into
a minimal topologically simple group. \qed
\label{embedding}
\end{corol}

\subsubsection*{Application to embeddings}
The last two results by Uspenskij provide a very interesting and
quite unexpected insight into a general problem about
embeddability of topological groups into groups with various additional
properties. 

The following question used to be quite popular
among the members of the school of topological algebra
headed by Alexander Arhangel'ski\u\i\ at Moscow University.
(Cf. \cite{Arh1,Arh2}.)
Suppose $\mathcal G$ is a certain class of topological groups.
Under which conditions is a given topological group $G$ isomorphic
with a topological subgroup of the direct product of a subfamily
of $\mathcal G$? In particular, is {\it every} topological group
$G$ isomorphic with such a subgroup? Using the notation adopted
in theory of varieties of topological groups \cite{Mvar}, the latter 
question can be restated as follows: given a class $\mathcal G$
of topological groups, when is it true that $SP({\mathcal G})$
contains all (Hausdorff) topological groups? 
(Here the letters $S$ and $P$ indicate the transition to a 
topological subgroup and the direct product, respectively.)

We believe that the following result (conjectured a few years
ago by the present author in \cite{Pbams}) appears here for the first time.

\begin{corol}
Let $\mathcal G$ be a class of topological groups such that
every topological group is isomorphic with a
topological subgroup of the direct product of a family of groups
from $\mathcal G$. Then every topological group is isomorphic
with a topological subgroup of a suitable group from $\mathcal G$.
\par
Put otherwise, if $SP({\mathcal G})$ is the class of all
topological groups, then already $S({\mathcal G})$ is the class of all
topological groups.
\end{corol}

\begin{proof}
Let $G$ be an arbitrary topological group. Without loss in
generality, assume that $G\neq\{e_G\}$.
Using Uspenskij's
Corollary \ref{embedding}, embed $G$ into a minimal, topologically
simple group $H$. Let now $H$ be embedded, as a topological subgroup,
into $\prod_{\alpha\in A}G_\alpha$, where $G_\alpha\in {\mathcal G}$.
For at least one $\alpha\in A$, the image of $H$ under the
$\alpha$-th coordinate projection $\pi_\alpha$ is non-trivial.
Because of topological simplicity of $H$, the kernel of
$\pi_\alpha\vert_H$ is $\{e_H\}$, and thus the restriction
$\pi_\alpha\vert_H$ is a (continuous) group monomorphism.
Because of minimality of $H$, the latter monomorphism is in fact
a topological isomorphism. Consequently, $\pi_\alpha\vert_G$
is a topological isomorphism of $G$ with a topological subgroup
of $G_\alpha$.
\end{proof}

\begin{corol}
If $\mathcal G$ is a class of topological groups closed under
formation of topological subgroups and having the property that
every topological group is isomorphic with a
topological subgroup of the direct product of a family of groups
from $\mathcal G$, then $\mathcal G$ is the class of all
topological groups. \qed
\end{corol}

The above corollaries from Uspenskij's results provide 
a rather efficient tool for handling
questions of the type mentioned at the beginning of the subsection.

\subsection{Representations in reflexive spaces}
Let us now re-examine Teleman's theorem \ref{teleman} again. 
It is evident from the proof that
the result has the following stronger form (which we have already
repeatedly used above).

\begin{thm}[Teleman,1957] 
\label{teleman2}
Every Hausdorff topological group $G$ embeds,
as a topological subgroup, into the group $\Iso(E)_s$ of
isometries of a suitable Banach space $E$ equipped with the
strong operator topology. \qed
\end{thm}

Can the above result be further strengthened? 

The first thing to be 
observed is that $E$ cannot be replaced with a {\it Hilbert} 
space. Indeed, there
are known to exist topological groups
possessing no nontrivial strongly continuous unitary representations. 

\begin{example}
The first example
of such kind seems to belong to Herer and Christensen
\cite{HC}, and the group in question is an abelian topological group.
Here we will convey the idea of the construction. 

A non-negative
real-valued function $\varphi$ defined on elements of a sigma-algebra
$\mathcal A$
of subsets of a set $X$ is called a {\it submeasure} if it is countably
subadditive (that is, always 
$\varphi(\cup_{i\in\N} A_i)\leq\sum_{i\in\N}\varphi(A_i)$),
monotone (that is, $\varphi(A)\leq\varphi(B)$ whenever
$A\subseteq B$), and satisfies $\varphi(\emptyset)=0$. 
A submeasure $\varphi$ 
is called {\it pathological} if it is not identically
zero and there exists no sigma-additive measure
on $(X,{\mathcal A})$ all of whose null sets are null sets 
with respect to $\varphi$. It can be shown that pathological submeasures
exist on every non-atomic sigma-algebra of sets. In particular,
there exists such a submeasure, $\varphi$, on the sigma-algebra of
Borel subsets of the Cantor set $X=\{0,1\}^\omega$, 
and moreover one can assume
the `regularity' condition: $\varphi(V)>0$ for every nonempty
open and closed subset of $X$. Equip the
vector space $C(X)$ with the topology of
convergence in submeasure $\varphi$: basic neighbourhoods of zero
are of the form 
\[V_{\e,\delta}:=
\{f\in C(X)\colon \varphi(\{x\in X\colon \vert f(x)\vert>\delta\})<\e\}, 
\e,\delta>0.\]
It is shown in \cite{HC} that
the additive topological group of the topological vector space
$C(X)$ with the above topology admits no nontrivial
strongly continuous unitary representations in Hilbert spaces.
\label{hererex}
\end{example}

\begin{example}
\label{banaex}
Later on, Banaszczyk \cite{Ba}
had shown the existence of abelian Banach--Lie groups without
nontrivial strongly continuous unitary representations. 
They are topological factor-groups of the
additive group of the separable Hilbert space $l_2$ by a suitably chosen
discrete subgroup $\Gamma$. 

Recall that a topological group $G$ together with a 
smooth structure modelled on a (possibly infinite-dimensional)
Banach space is called a {\it Banach--Lie group}. 
The above smooth structure is determined by a neighbourhood of
identity, $V$, in $G$ and a homeomorphism $\phi$ between $V$ and
the open unit ball $B$ in a suitable Banach space $E$ satisfying
a certain smoothness condition. To describe it, 
notice that for each $g\in G$ the
formula $_g\phi(h):=\phi(g^{-1}h)$ determines a homeomorphism
between the neighbourhood $gV$ of $g$ and the ball $B$. The domain of
the composition map $\phi\circ(_g\phi)^{-1}$ is therefore an open
subset of $V$ (possibly empty). The smoothness condition now
is this: whenever $g\in G$, the map $\phi\circ(_g\phi)^{-1}$ is
$C^\infty$ in its domain of definition.
In particular, it follows easily that a topological factor group of
a Banach-Lie group by a discrete (though not necessarily by an arbitrary
closed) subgroup inherits the Banach-Lie
structure in a canonical way.
A good introduction to the theory of Banach--Lie groups along these
lines can be found in Karl Hofmann's lecture notes \cite{Hof},
and some of it made its way into the recently published
book \cite{HM}.

Banach-Lie groups form in a sense the closest class of topological
groups (with additional structure) to that of finite-dimensional
Lie groups, and such a drastic difference in behaviour between
groups from two classes is rather stunning. 
Indeed, it is worth recalling a well-known and easy to prove fact:
every locally compact group
admits a faithful strongly continuous unitary representation. Such
is the left regular representation of $G$ in the space $L_2(G)$
of all square-integrable complex-valued functions on $G$ with
respect to the Haar measure. Locally compact abelian (LCA) groups
enjoy a much stronger property: characters (that is, continuous
homomorphisms to the circle rotation group ${\mathbb T}=U(1)$)
separate points in LCA groups. Notice that the groups in Banaszczyk's
example are abelian, and having no nontrivial unitary representations
is a much more restrictive property than just having no
continuous characters.
\end{example}

The next natural thing to ask is, can one at least assume that
the Banach space $E$ in Theorem \ref{teleman2} 
is reflexive? 

A negative answer was very recently given by Megrelishvili 
\cite{Meg3} who has
shown that the topological group $\Homeo_+({\mathbb I})$
of all orientation-preserving homeomorphisms of the
closed interval equipped with the
compact-open topology admits no non-trivial representations in 
reflexive Banach spaces. The rest of this entire section is
loosely grouped around a sketch of the proof of this
result. 

We will begin with the following criterion. Recall
that the weak topology on a topological vector space $E$ is
the coarsest topology with respect to which 
every continuous functional $\phi$ on
$E$ remains continuous. The space $E$ equipped with the weak
topology will be denoted by $E_w$. Recall also
that a (real or complex-valued) bounded
function $f$ on a topological group $G$ is called {\it weakly almost
periodic} ({\it WAP}) \cite{Eb}, 
if the set of all left translations
$\{L_gf\colon g\in G\}$ of this function is weakly relatively compact
in the space $C^b(G)$ of all continuous bounded functions on
$G$ with the topology of uniform convergence. In other words, 
the closure of $\{L_gf\colon g\in G\}$ in the space $C^b(G)_w$ is
compact. It is useful to notice that the concept of
weak almost periodicity is partially 
independent of the topology on the group
$G$ in the sense that the same function $f$ remains 
weakly almost periodic on
the group $G$ equipped with the discrete topology.

\begin{thm}[Shtern, 1994, \cite{Shtern} and Megrelishvili, 1998,
\cite{Meg2}] 
A topological group $G$ embeds into the isometry
group of a reflexive Banach space equipped with the strong operator
topology if and only if weakly almost periodic functions
separate points and closed subsets in $G$.
\label{wap}
\end{thm}

The proof will be preceded by a few concepts. Suppose a group $G$ acts
on a normed space $E$ by isometries, that is, we are given a
homomorphism $\pi\colon G\to\Iso(E)$ (a representation of $G$). 
Fix a vector $\xi\in E$ (usually
of norm one) and a bounded linear functional $\phi$ on $E$.
The function 
\[\tau_{\xi,\phi}:G\ni g\mapsto \phi(g\cdot\xi)\in\K\]
(where $\K=\R$ or $\C$)
is called the $(\xi,\phi)${\it -th matrix coefficient} (of the
representation $\pi$). The topology induced on $G$ 
(or else on the group of isometries $\Iso(E)$) by the family of
all matrix coefficients is called the {\it weak operator topology}.
If $E$ is reflexive, then the weak operator topology and the strong
operator topology on $G$ coincide. This important fact was known for
the Hilbert spaces $E$ since long ago, 
but for general reflexive Banach spaces it
was only recently established by Megrelishvili \cite{Meg2}.
(Notice that for non-reflexive spaces there are counterexamples
distinguishing between the two topologies on groups of isometries,
{\it loco citato.})

It is easy to notice that if $E$ is a reflexive Banach space,
then every matrix coefficient $\tau_{\xi,\phi}$
is a weakly almost periodic function on $G$. 
Firstly notice that the $G$-orbit of this coefficient
consists of all matrix coefficients of the form 
$\tau_{\xi,g\cdot\phi}$, $g\in G$. Further, the mapping
\[E^\ast\ni\chi\mapsto \tau_{\xi,\chi}\in C^b(G)\]
is a bounded (= continuous) linear operator 
with respect to the norm topologies
on both spaces (here $\xi\in E$ is fixed), therefore it remains
continuous relative to the weak topologies on both spaces in
question and, in particular, the restriction
\[B^\ast_w\ni\chi\mapsto \tau_{\xi,\chi}\in C^b(G)_w\]
is continuous, where both the unit ball in the dual space and
the space of continuous functions are equipped with their weak
topologies. The weak closure of the $G$-orbit of a unit functional
$\phi$ is compact (the reflexivity of a Banach space is equivalent
to the weak compactness of the dual unit ball). Consequently,
the image of this closure in $C^b(G)_w$ is compact. But this is
exactly the weak closure of
the $G$-orbit of the matrix coefficient in question.

Hence follows the {\it necessity} in Theorem \ref{wap}: suppose
$E$ is a reflexive Banach space and $G\hookrightarrow\Iso(E)_s$
is a topological embedding. According to the above result by Megrelishvili,
matrix coefficients generate the topology of $G$. It remains
to notice that the collection $\mathrm{WAP}(G)$ of all weakly almost periodic
functions on $G$ (which in fact forms a $C^\ast$-algebra) is closed under
the lattice operations, such as taking maxima of finite families
of functions. The rest is just simple routine of general topology.

Now let us outline the proof of {\it sufficiency} of Th. \ref{wap}
as proposed by
Megrelishvili \cite{Meg2}. (The proof of Shtern drew upon
algebra representation theory.) 
Make the collection $\mathrm{WAP}(G)$ of all weakly almost periodic functions 
on a topological
group $G$ into a Banach space by equipping it with the supremum norm. 
Let $f\in \mathrm{WAP}(G)$ be arbitrary, and denote by
$E=E_{f}$ the linear span of  
the orbit $G\cdot f$ in $\mathrm{WAP}(G)$. Denote by $W=W_f$ the 
convex circled envelope of $G\cdot f$ in $E$. (Recall
that a subset $A$ of a vector space is {\it circled} if
$\lambda A\subseteq A$ whenever $\lambda$ is a scalar with
$\abs\lambda\leq 1$; for real vector spaces circled sets are
just symmetric sets.)
As a consequence of
the Krein--Smulian Theorem, $W$ is relatively weakly compact
(indeed, for every $f\in \mathrm{WAP}(G)$, the orbit
$G\cdot f$ is relatively weakly compact).

Now a procedure developed in \cite{DFJP} is applied.
For every natural $n$, let $\norm\cdot_n$ denote the equivalent 
norm on the space $E$ whose
open unit ball is
\[U_n:= 2^nW+2^{-n}B_E,\]
where $B_E$ is the unit ball of $E$ with respect to the induced
(supremum) norm. Define
\[\norm x_\ast:=\left(\sum_{n=1}^\infty \norm x_n^2\right)^{\frac 12}\]
for every $x$ from the set $X=X_f$ of elements for which the expression
above assumes a finite value. The norm $\norm\cdot_\ast$
is invariant under the action of $G$ by left translations
(as $X$ is a translation-invariant collection of 
functions on $G$). The weak compactness of $W$ results in the weak
compactness of the unit ball of $\norm\cdot_\ast$, which fact
implies that
$X$ is a reflexive Banach space. \par
The space $X$ is 
continuously embedded into $E$ in a canonical way so
that on every bounded subset of $X$ the embedding is
a homeomorphism with respect to the weak topologies.
Since the orbit $G\cdot f$ is bounded in $X$, it follows that
the representation of $G$ in $X$  is a strongly continuous
representation by isometries. \par
Form the $l_2$-direct sum
of all representations of $G$
obtained in the above way as $f$ runs over
the collection of all weakly almost periodic functions on $G$.
One obtains a representation of $G$ in a reflexive Banach space 
$Y=\oplus_{f\in \mathrm{WAP}(G)}X_f$ by isometries which
determines a topological group embedding
$G\hookrightarrow\Iso(Y)_s$.

\begin{remark}
Of course the above criterion \ref{wap}
can be slightly adjusted to suit a variety of different
situations. For example, a topological group $G$ admits a separating
family of strongly continuous representations in reflexive Banach
spaces if and only if WAP functions separate points in $G$.
Or else: the existence of a non-trivial representation in a reflexive
Banach space is
equivalent to the existence of a non-constant WAP function on $G$.
\end{remark}

\subsection{\label{compacti}Compactifications of topological groups}
Now we can proceed to an example of a topological group admitting
no non-trivial strongly continuous representations by isometries in
reflexive Banach spaces. It is interesting to note that such an
example is supplied by one of the most common `massive' groups.

\begin{thm}[Megrelishvili, 1999, \cite{Meg3}] 
Every weakly almost periodic function on the topological group
$\Homeo_+[0,1]$ is constant. Consequently, this group admits no
non-trivial strongly continuous representations in reflexive
Banach spaces.
\label{homeo}
\end{thm}

To better understand the idea of the proof, it is useful to start with
the concept of the {\it lower uniformity} on a topological group --
a concept that has growing significance on its own, especially, it seems,
for `massive' topological groups, and which was first investigated
by Roelcke (see e.g. \cite{RD}) and since then
brought to the limelight through the work of Uspenskij. 
\par Let us begin with
an obvious observation: both standard uniformities on
a topological group $G$ --- namely, the right uniformly equicontinuous
uniform structure ${\mathcal U}_\Rsh(G)$ and its left counterpart
${\mathcal U}_\Lsh(G)$ --- are {\it compatible uniformities}, that is, 
each of them generates the
topology of $G$. As a straightforward consequence,
the {\it upper} uniform structure,
which is the supremum of the two,
\[{\mathcal U}_\vee(G):={\mathcal U}_\Rsh(G)\vee {\mathcal U}_\Lsh(G),\]
is compatible as well.
The upper uniform structure is best known in the context of
{\it completeness} of topological groups. Indeed, not every
topological group embeds into a topological group that is complete
with respect to any (equivalently: both) 
of the one-sided structures ${\mathcal U}_\Rsh(G)$
and ${\mathcal U}_\Rsh(G)$ (such groups are called 
{\it Weil-complete}.) A counter-example 
(Dieudonn\'e \cite{Dieu}) is the same group $\Homeo_+[0,1]$.
At the same time, every topological group embeds as an everywhere
dense subgroup into a group that is complete relative to the upper
uniformity. (Which fact of course just means that the concepts of
completeness and completion with respect to one-sided uniformity
are misfits that must be discarded in favour of
two-sided completeness and completion.)

The {\it lower} uniformity is the infimum of the left and right
equicontinuous uniformities:
\[{\mathcal U}_\wedge(G):={\mathcal U}_\Rsh(G)\wedge {\mathcal U}_\Lsh(G).\]
It is also a compatible uniformity, which fact
can be seen from the explicit form of basic entourages,
\[V_\wedge:= \{(x,y)\in G\times G\colon x\in VyV\},\]
where $V$ runs over neighbourhoods of the identity element $e_G$.

A remarkable observation is that for many `massive'
topological groups the lower uniform structure 
${\mathcal U}_\wedge(G)$ turns out to be precompact. 
Uspenskij calls topological groups with this property
{\it Roelcke-precompact.} We side with the author of the review
\cite{Grant} and adopt a more
functional terminology, calling such groups
{\it lower precompact,} and their completions with respect to the
lower uniformity {\it lower completions.}

Here are some of the major examples of lower precompact topological
groups.

\begin{example} The full unitary group $U(\H)_s$ of an
infinite-dimensional Hilbert space with the strong operator
topology. The lower compactification of $U(\H)_s$ can be
identified with the semigroup of all operators on $\H$ of
norm $\leq 1$, equipped with the weak operator topology. \cite{U4}.
\end{example}

\begin{example} Let $X$ be an infinite set. 
Equip the full group of permutations $S(X)$ 
with the topology of simple convergence (where $X$ is 
regarded as a discrete set), that is, the topology induced
from the Tychonoff power $X^X$. Then $S(X)$ is lower-precompact.
\cite{RD}
\end{example}

\begin{example}
The group $\Homeo_+(I)$ of orientation-preserving homeomorphisms
of the closed unit interval with the compact-open topology.
Identify each homeomorphism $h\colon I\to I$ with its graph in
the square $I\times I$.
Then the lower completion of the group $\Homeo_+(I)$
can be identified with the collection of all $C^0$ curves
$\gamma$ in $I\times I$ starting at the lower left corner $(0,0)$
and ending at the right upper corner $(1,1)$ and never going
either to the left or down --- more exactly, if
an orientation-preserving parametrization of $\gamma$
is chosen and
two values of the parameter satisfy $t_1\leq t_2$, one necessarily has
$\gamma(t_1)_i\leq\gamma(t_2)_i$, $i=1,2$. (Fig. 1.)
\bigskip

\begin{center}
\scalebox{0.3}[0.3]{\includegraphics{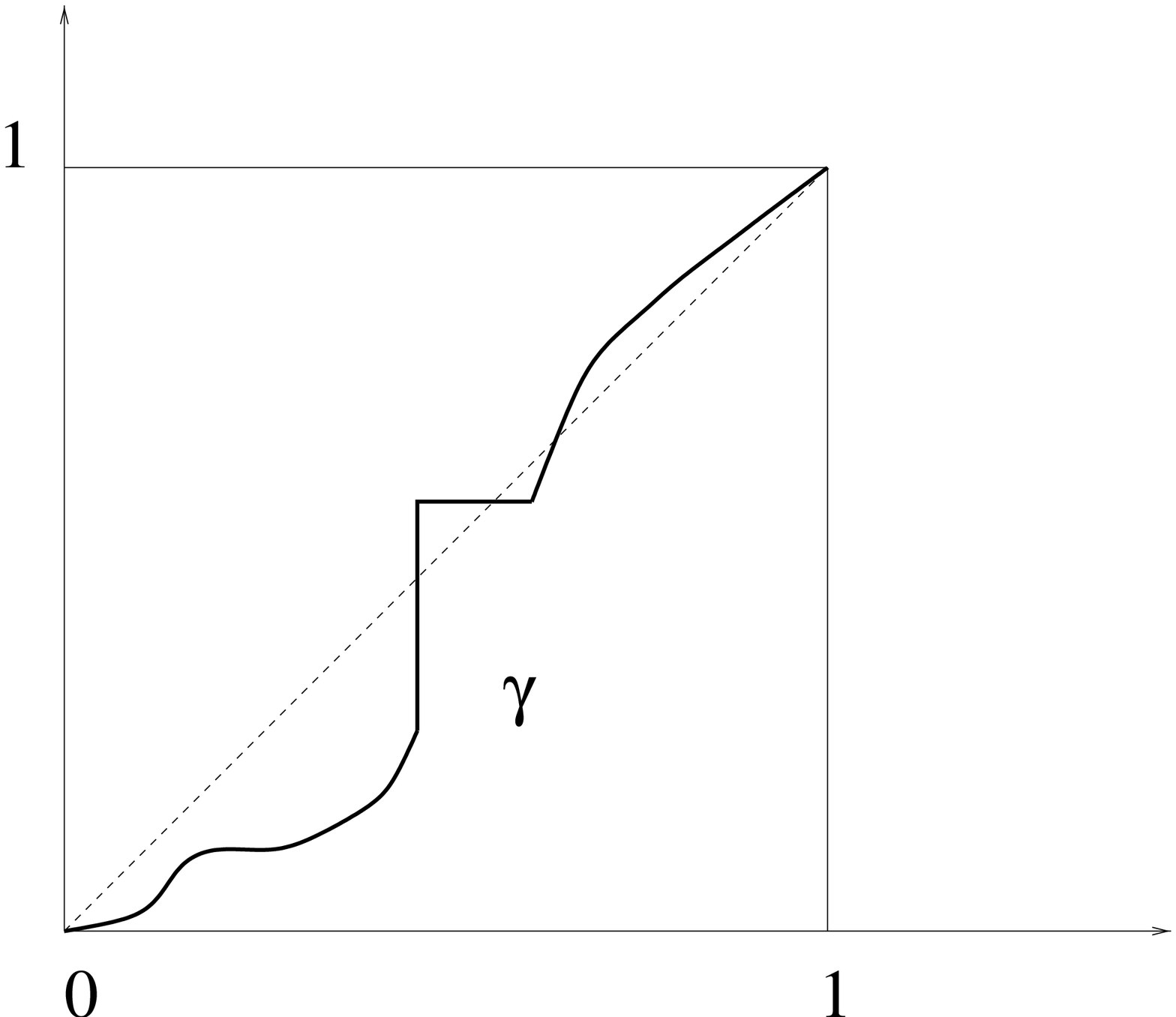}}\\
{\small
Fig. 1.
An example of an element of $\gamma_\wedge(\Homeo_+({\mathbb I}))$.}
\end{center}
\bigskip
The topology on the
set of curves is that of uniform convergence.
(Uspenskij, \cite{U5}.)
\end{example}

\begin{example} Let $M$ be a complete Urysohn space that is
$\omega${\it -homogeneous}, that is, every isometry between two
finite metric subspaces of $M$ extends to an isometry of $M$
onto itself. (In the separable case, $\omega$-homogeneity of a complete
Urysohn metric space can be taken for granted.)
The full group of isometries $\Iso(M)_s$ of $M$ 
equipped with the pointwise topology is lower precompact.
(Uspenskij, \cite{U3}.) In this case, the lower completion of
$\Iso(M)_s$ can be identified with a set of functions
$M\times M\to I$ satisfying a certain set of conditions which we
do not reproduce here.
\end{example}

The lower completion of a lower precompact topological group $G$ always
supports the structure of
a $G$-space, and sometimes that of a semigroup, 
which structures can be used to 
establish certain topological-algebraic
properties of $G$ itself such as minimality and topological
simplicity. We invite the interested reader to consult 
Uspenskij's papers \cite{U4,U3}.

\begin{remark}
There is an interesting yet largely unexplored connection between
the above concepts of the lower precompactness and lower completion,
on the one side,
and modern representation theory of `massive' groups, on the other.
It appears that many major examples of infinite
dimensional groups admitting a tractable representation
theory in Hilbert spaces and commanding prominence in
mathematical physics are lower precompact if equipped with
a `natural' topology. 

Namely, in a theory developed over the past decade or two,
largely through the efforts of Neretin, G. Ol'shansky and others
(the book \cite{Ner} is the most up-to-date source), to every
`massive' group there is associated a certain compact semigroup
called the {\it mantle} of $G$ and having the same representation
theory in Hilbert spaces as $G$. The mantle of a group happens
to coincide with the lower
completion in the most telling case of 
the full unitary group. One certainly expects
the lower completion of a lower precompact group
to be typically `larger' than the mantle,
as the group $\Homeo_+(I)$ examplifies (whose mantle is of course
trivial, as one of the consequences of Megrelishvili's result).
Yet it may well happen
that for many important lower precompact groups the mantle 
coincides with the lower completion.

Among the topological groups whose mantle has been computed
are groups of
diffeomorphisms of manifolds, groups associated to Virasoro
and other Kac--Moody algebras, various
groups of operators in Hilbert spaces (subgroups of unitary
groups preserving one 
structure or other), groups of currents, and groups of
automorphisms of measure spaces. It is interesting to examine
each of these groups for lower precompactness
and to compare the lower completion with
the mantle.
\end{remark}

What is the relationship between lower precompactness and
weak almost periodicity?
It turns out that every weakly almost periodic
function $f$ on a topological group $G$ is uniformly
continuous with respect to both uniformities ${\mathcal U}_\Rsh(G)$,
${\mathcal U}_\Lsh(G)$. (For a proof of this fact, see e.g. 
\cite{Ruppert}. Notice that one half of this statement is
nearly obvious: the ${\mathcal U}_\Lsh(G)$-uniform
continuity follows in a straightforward fashion from the
compactness of the orbit of $f$ in the {\it pointwise} topology.)
Consequently, the algebra
$\mathrm{WAP}(G)$ is contained in the algebra $C^b_\wedge(G)$ of all
bounded lower uniformly continuous functions.

Now recall that every algebra $A$ formed by bounded continuous functions on a
topological space $X$ determines a compactification of $X$, which is
the completion with respect to the coarsest uniform structure
making each function $f\in A$ uniformly continuous.
We will denote this compactification by $\gamma_A(X)$ (though in fact
$\Spec A$ would be more informative and commonly recognizable).
The compact space $\gamma_A(X)$ is the maximal ideal space of $A$.
Teleman's argument actually tells us that in the case where 
$X=G$ is a topological group, the algebra $A$ is
closed under translations by elements of $G$, and the action of $G$
on $A$ by left translations is continuous (where $A$ is equipped with
the topology of uniform convergence on $G$), then the resulting action
of $G$ on the compactification $\gamma_A(G)$ is continuous as well,
that is, $\gamma_A(G)$ is a compact $G$-space. (Indeed, 
one only has to notice that in this case $A$ is contained in the
algebra $C^b_\Rsh(G)$ as a $G$-invariant subalgebra, and
the compactification $\gamma_A(G)$ in question forms
a $G$-invariant compact subset of the unit ball $B_w^\ast$ in the dual
space $C^b_\Rsh(G)'$; since the action of $G$ on the ball is continuous,
so is the action of $G$ on $\gamma_A(G)$.)

Let us consider some examples of topological group
compactifications of this sort.

\begin{example}
Here is the master concept.
A bounded scalar-valued function $f$ on a 
topological group $G$ is called
{\it almost periodic} if the $G$-orbit of $f$ is relatively compact in
$C^b(G)$. The collection of all almost periodic functions is denoted by
$\mathrm{AP}(G)$, and the corresponding compactification 
$\gamma_{\mathrm{AP}}(G)$ is
known as
the {\it Bohr compactification} of $G$. The Bohr compactification of a
topological group is itself
a compact group, and it is maximal among
all compact groups into which $G$ admits a continuous homomorphism with
a dense image. This concept is one of the cornerstones of the classical
abstract harmonic analysis \cite{HR}. 
The following two related notions will be
used later on. A topological group $G$ is {\it maximally almost
periodic} ({\it MAP})  if the canonical continuous
homomorphism $G\to\gamma_{\mathrm{AP}}(G)$ is a
monomorphism, and {\it minimally almost periodic} ({\it map}) if 
$\gamma_{\mathrm{AP}}(G)=\{e\}$. These can be restated in terms
of representation theory as follows:
$G$ is MAP iff continuous
unitary representations in finite-dimensional
Hilbert spaces separate points in $G$, and $G$ is map iff it possesses
no non-trivial continuous finite-dimensional unitary representations.
If $G$ is abelian, then `finite-dimensional unitary representations' in
the above criteria can be replaced with `characters' (continuous
homomorphisms to the circle rotation group ${\mathbb T} = U(1)$).
\end{example}

\begin{example} Of singular importance in abstract topological dynamics
\cite{Aus, dV1, P19}
is the compactification $\gamma_\Rsh(G)$ formed
with respect to the algebra
$A=C^b_\Rsh(G)$. This compactification is known as the
{\it greatest ambit} of $G$ and possesses a certain universal property
which we will now establish. Let $G$ be a topological group. A compact
$G$-space $X$ together with a distinguished point $x^\ast\in X$ is
called an {\it ambit} if the orbit $G\cdot x^\ast$ is everywhere
dense in $X$. It is clear how to define morphisms between two
$G$-ambits $(X,G,\tau_X,x^\ast)$ and $(Y,G,\tau_Y,y^\ast)$:
such a morphism is a continuous map $\varphi\colon X\to Y$ which is
equivariant, that is, commutes with the actions $\tau_X$ and $\tau_Y$
in the sense that, for all $x\in X$ and $g\in G$,
\[\varphi(g\cdot_X x) = g\cdot_Y \varphi(x),\]
and also $\varphi$ preserves the distinguished points:
\[f(x^\ast)=y^\ast.\]
One can make the compactification $\gamma_\Rsh(G)$ into an ambit
by marking as the distinguished element $e_G^\ast$, which is
simply the identity element of $G$, or rather its image in the
compactification. 
Since the embedding
$G\hookrightarrow \gamma_\Rsh(G)$ is clearly topological 
(the bounded ${\mathcal U}_\Rsh(G)$-uniformly continuous functions
separate points and closed subsets in $G$), we can identify $G$ with
a topological subspace of the greatest ambit. 
For every $G$-ambit $(X,x^\ast)$ there is a unique
morphism of $G$-ambits $\varphi\colon \gamma_\Rsh(G)\to X$. 
Indeed, for every continuous (real or complex-valued)
function $f$ on such a $G$-ambit $X$, the pullback $f_\ast$ of
$F$ to $G$ defined for each $g\in G$ by
\[f_\ast(g):=f(g\cdot x^\ast)\]
is ${\mathcal U}_\Rsh(G)$-uniformly continuous (an easy, direct check).
Moreover, the map $f\mapsto f_\ast$ as above is a homomorphism of
$C^\ast$-algebras from $C(X)$ to $C^b_\Rsh(G)$ (that is, it is linear,
bounded, multiplicative, and -- in the complex case -- preserves
the involution). This map determines a continuous mapping between the
corresponding compactifications going in the opposite direction, which
is exactly 
the morphism of $G$-ambits $\gamma_\Rsh(G)\to X$ we are after.
Often the greatest ambit of a topological group
$G$ is denoted by ${\mathcal S}(G)$. 

While of course there is a canonical morphism of $G$-ambits
from $\gamma_\Rsh(G)$ onto the Bohr compactification $\gamma_{\mathrm{AP}}(G)$,
it is one-to-one if and only if $G$ is a precompact group.
\end{example}

\begin{example}
The compactification $\gamma_{\mathrm{WAP}}(G)$ of a topological group with
respect to the algebra $\mathrm{WAP}(G)$ of weakly almost periodic functions
is known as the {\it maximal 
semitopological semigroup compactification of $G$} and possesses a
universal property similar to that of the greatest ambit with
respect to all compact semitopological (that is, the operations
are {\it separately} continuous) semigroups containing an image of
$G$ as an everywhere dense subsemigroup. Clearly, there is a
canonical $G$-equivariant mapping onto 
$\gamma_\Rsh(G)\to\gamma_{\mathrm{WAP}}(G)$. Megrelishvili's Theorem
\ref{homeo} says, equivalently, that the compactification
$\gamma_{\mathrm{WAP}}(\Homeo_+{\mathbb I})$ is a singleton!
\end{example}

\begin{example}
If $A=C^b_\wedge(G)$ is the algebra of all bounded lower
uniformly continuous functions on a topological
group $G$, then we will denote the corresponding
compactification of $G$ by $\gamma_\wedge(G)$.
A topological group $G$ is lower precompact if and only if the
compactification mapping $G\to\gamma_\wedge(G)$ induces on $G$ the
lower uniformity. There is a canonical 
continuous mapping from 
$\gamma_\wedge(G)$ onto the compactification $\gamma_{\mathrm{WAP}}(G)$.
Every weakly almost periodic function $f$ uniquely extends from $G$ to
$\gamma_{\mathrm{WAP}}(G)$, and therefore clearly factors through the mapping
$\gamma_\wedge(G)\to \gamma_{\mathrm{WAP}}(G)$. In other words, weakly almost
periodic functions on $G$ can be identified with those continuous
functions on the lower compactification $\gamma_\wedge(G)$ whose orbit under
left translations by elements of $G$ 
is weakly relatively compact
in the Banach space $C(\gamma_\wedge(G))\cong C^b_\wedge(G)$. 
For those lower precompact groups whose lower completion 
admits a transparent geometric interpretation, this observation makes
working with weakly almost periodic functions (notoriously evasive
objects) somewhat easier. In particular, this applies to the group
$\Homeo_+({\mathbb I})$. 
\end{example}

\begin{remark} A reference to the subject of compactifications of
the kind described above is the book \cite{BJM}.
\end{remark}

Now we are in a position to convey the flavour of
Megrelishvili's proof of Theorem \ref{homeo}. The proof is remarkably
`graphical.' For each 
triple of real numbers $a,b,c$ with $0\leq a\leq c \leq b\leq 1$
denote by $\beta^{a,c,b}$ the element of the lower completion
$\gamma_\wedge(\Homeo_+(\I))$ represented by the following curve
(Fig. 2).
\bigskip

\begin{center}
\scalebox{0.3}[0.3]{\includegraphics{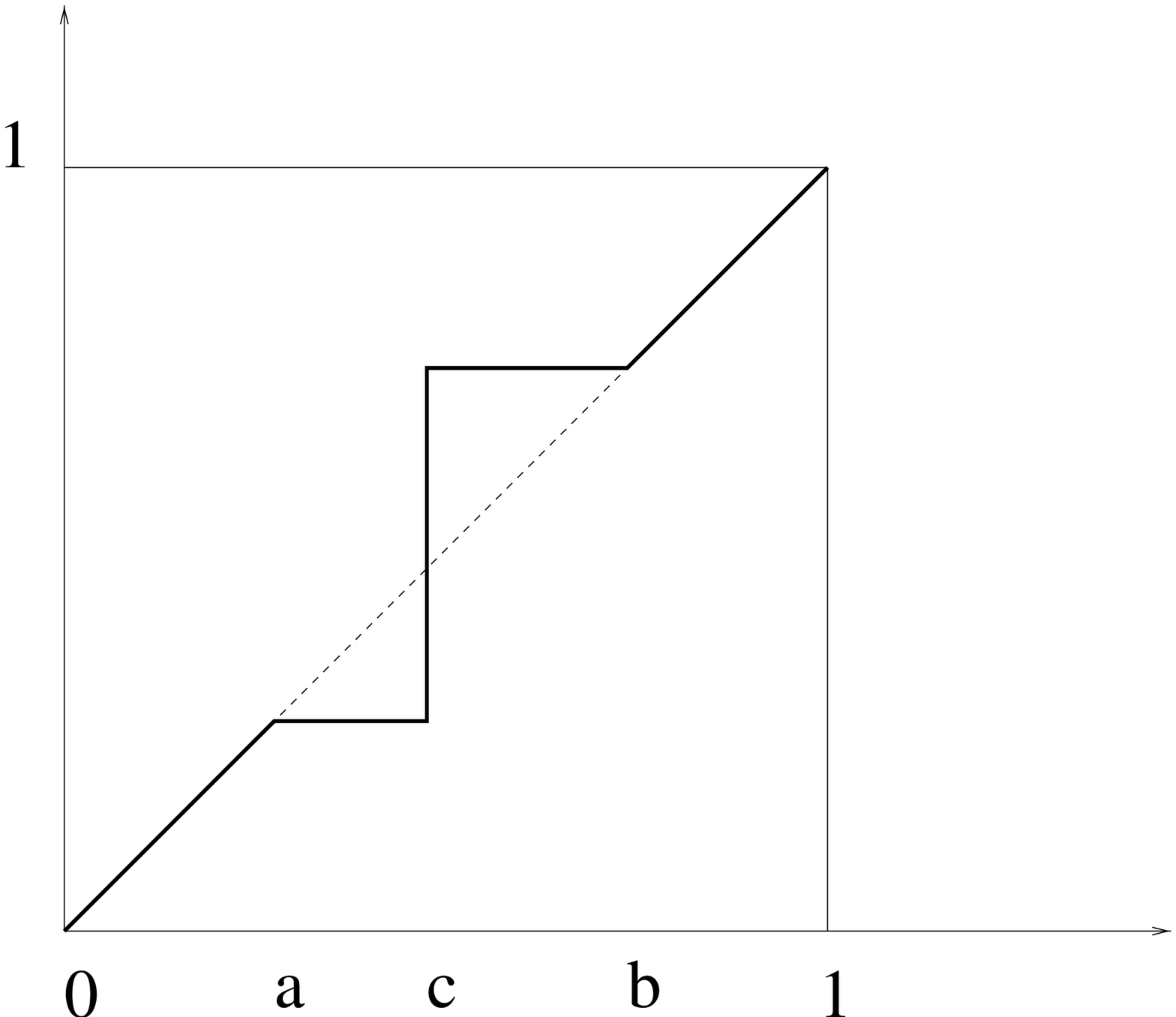}}\\
{\small
Fig. 2.
The curve $\beta^{a,c,b}$.}
\end{center}
\bigskip
As we have already remarked, the compactification
$\gamma_{\mathrm{WAP}}(\Homeo_+(\I))$ supports a canonical semigroup
structure. It turns out that 
the image of the curve $\beta^{0,1,1}$ in
$\gamma_{\mathrm{WAP}}(\Homeo_+(\I))$ forms the zero element.
More generally, using the fact that the latter semigroup is
semitopological, one can verify that {\it every} zig-zag curve 
made up entirely of alternating horizontal and vertical segments, that is,
one of the form 
\[\beta^{a_0,a_1,a_2}\circ \beta^{a_2,a_3,a_4}\circ\cdots
\beta^{a_{n-2},a_{n-1},a_n},\] 
where $a_0,a_1,\cdots,a_n$ is a partition of the interval $\I$, 
goes to zero element under the map 
\[\gamma_\wedge(\Homeo_+(\I))\to\gamma_{\mathrm{WAP}}(\Homeo_+(\I)).\]
At the same time, as the mesh of the partition tends to zero, the 
zig-zag curves as above
converge to the identity homeomorphism and consequently
their images in the semigroup $\gamma_{\mathrm{WAP}}(\Homeo_+(\I))$ 
must converge to the
identity element. But every semigroup whose zero and identity
elements coincide
with each other is trivial. Q.E.D.

\begin{remark}
We hope that the above hands-waving argument has not left an
impression of Megrelishvili's proof being easy: a precise incarnation
of the idea happens to be very involved technically.
\end{remark}

\begin{question}
(Megrelishvili; Akin and Glasner.) 
Does there exist a {\it monothetic} 
topological group whose points and
closed subsets are not separated by WAP functions?
\end{question}

As observed by Megrelishvili, this is equivalent to the same question
about arbitrary 
{\it abelian} topological groups, in view of Theorem \ref{ab} below.

We conclude this Section by recalling a long-standing open problem.

\begin{problem} (Shtern, \cite{Shtern0}) 
Give an intrinsic description of those
topological groups $G$ admitting a separating family of strongly
continuous unitary representations in Hilbert spaces (or else
embeddable, as topological subgroups, into the unitary groups of
Hilbert spaces with the strong operator topology).
\label{questshtern}
\end{problem}

It is worthwhile looking at a well-known description of such
topological groups --- which however falls short of an acceptable
intrinsic criterion --- after which the very Theorem \ref{wap} 
was fashioned. 
Recall that a complex-valued function $f$ on a
group $G$ is called {\it positive definite} if for every finite
collection $g_1,g_2,\cdots,g_n$ of elements of $G$ and every
collection of the same size of complex numbers, $\lambda_1,\lambda_2,
\cdots,\lambda_n$, one has
\[\sum_{i,j=1}^n \lambda_i\overline{\lambda_j}f(g_i^{-1}g_j)\geq 0.\]
The relation between positive definite functions and unitary
representations is as follows. If $\pi\colon G\to U(\H)$ is a
unitary representation of a group, then for every $\xi\in\H$ the
function 
\[G\ni g\mapsto (g\cdot\xi,\xi)\in\C\]
is positive definite. (Of course, this function is nothing but
the matrix coefficient $\tau_{\xi,\hat\xi}$, where $\hat\xi$
denotes the linear functional represented by $\xi$ in the Hilbert
space.)
Conversely, given a positive definite function
$f$ on $G$, one can construct a unitary representation of $G$ as
follows: the Hilbert space $\H$ is the completion of the vector space
$\mathrm{lin}\,(G)$
spanned by elements of $G$ as a Hamel basis and equipped with the
inner product defined on elements of $G$ by the formula
$(g,h):= f(g^{-1}h)$ and extended all over $\mathrm{lin}\,(G)$ by
sesquilinearity. (Positive definiteness of $f$ serves to verify
the first axiom of the inner product.)

Using the above two constructions, is not difficult to prove that a 
topological group $G$ embeds into the unitary
group $U(\H)_s$ of a Hilbert space if and only if 
continuous positive definite functions on $G$ separate points and elements
of some closed subbase for $X$.
Even if this description (along with obvious variations)
is being extensively used in functional analysis and representation
theory, it clearly falls short of what one would accept as an answer
to Shtern's question \ref{questshtern}.
The theorem of F\o lner--Cotlar--Ricabarra 
(Th. \ref{fcrek} below) expemplifies
what is being accepted as a satisfactory answer to questions of the
above kind.

It is interesting that the representability in reflexive Banach
spaces and unitary representability of topological groups have
not been distinguished from each other either.

\begin{question}[Shtern, \cite{Shtern}] Is it true that a topological
group $G$ admits a compete system of strongly continuous 
representations by isometries in reflexive Banach spaces if and
only if $G$ admits a complete system of strongly continuous
unitary representations?
\end{question}

Two obvious candidates for counter-examples are groups constructed
by Herer--Christensen (Ex. \ref{hererex}) 
and by Banaszczyk (Ex. \ref{banaex}): how many WAP functions
do they possess?

\section{Free actions {\it vs} fixed point on compacta property}

\subsection{Veech's theorem}
Now that we have examined the question of existence of effective
actions on compacta along with a range of broader topics prompted by
this question, it is natural to ask:
what about the existence of {\it free} actions?

Rather remarkably, an answer to this question draws a watershed between the
locally compact groups and many concrete
`massive' topological groups.
By far the most important 
result in the affirmative direction is the following.

\begin{thm}[W. Veech, 1977, \cite{V1}] 
Every locally compact group $G$ acts freely on a compact space.
\end{thm}

The proof we shall now outline is based on a simplification of Veech's
original argument proposed by Pym \cite{Pym}.
First of all, it follows from the 
`general nonsense' of the theory of uniform
compactifications that if $Y\subseteq (G,{\mathcal U}_\Rsh(G))$ 
is a uniformly discrete subspace, then
the closure of $Y$ in the greatest ambit $\gamma_\Rsh(G)$ is
canonically homeomorphic to the Stone-\u Cech compactification
$\beta X$ of the space $X$ with the discrete topology. 
The property of a subset $X$ of a topological group $G$ being
${\mathcal U}_\Rsh(G)$-uniformly discrete means that for some
neighbourhood $U$ of the identity:
\begin{equation}
Ux\cap Uy=\emptyset \mbox{ for all
$x,y\in X$, $x\neq y$.}
\label{discrete}
\end{equation}
Let $V$ be an open set with $\mathrm{cl}\,V$ compact
and contained in $U$.
It follows that the open set $V\cdot\mathrm{cl}\,X$ in
$\gamma_\Rsh(G)$ is canonically homeomorphic to the product 
$V\times\beta X$.

Given a symmetric compact
neighbourhood $U\ni e_G$, one can
construct a maximal subset $X$ with the property
(\ref{discrete}).
Clearly, the sets $U^2x$, $x\in X$ form a cover of $G$, and
consequently $\mathrm{cl}\,U^2\cdot\mathrm{cl}\,X$ coincides with
all of $\gamma_\Rsh(G)$. (Notice: here the compactness of 
$\mathrm{cl}\,U^2$ is used in an essential sort of way, and
this is precisely where the argument breaks down 
for more general topological groups.)
If now $x^\ast\in\gamma_\Rsh(G)$ is
arbitrary, then $x^\ast=u\cdot x$ for some 
$u\in\mathrm{cl}\,U^2\subseteq G$
and $x\in\mathrm{cl}\,X$, and therefore $x^\ast$ belongs to the
closure of the set $uX$ which is a maximal uniformly discrete set
with respect to the neighbourhood of identity $uUu^{-1}$.

Let $g\in G$, $g\neq e_G$, and let
$x^\ast\in\gamma_\Rsh(G)$ be arbitrary.
As we have seen, there is no loss in
generality is assuming that $x^\ast$ belongs to the closure of some
set $X\subseteq G$ which is maximal with respect to the 
property (\ref{discrete}), where $U$ is a suitable compact neighbourhood
of identity. Of course one can choose $U$ as small as desired, in
particular satisfying $g\notin U^2$.

Let $e_G\in V^2\subset U$.
Elementary combinatorial considerations 
show that $X$ can be
partitioned into finitely many pieces $X_1,\dots,X_k$ in such
a way that
for each $i=1,2,\dots,k$ the sets $V\cdot X_i$ and $g\cdot X_i$ are
disjoint in $G$. 
From the above description of the  
topological structure of $U\cdot \mathrm{cl}\,X$ in $\gamma_\Rsh(G)$
it follows that the set 
$\mathrm{cl}\,X_i$ is disjoint from its translation by $g$.
Since $x^\ast\in\mathrm{cl}\,X_{i_0}$ for some $i_0$, we conclude
that $g\cdot x^\ast\neq x^\ast$. Q.E.D.

Recall that an action of a group $G$ on a topological space $X$
is called {\it minimal} if the orbit of every point is
everywhere dense in $X$.
It is easy to see that 
if a topological group $G$ admits a free action on a compact space,
then $G$ admits a free (and {\it ipso facto} effective) 
{\it minimal} action on a compact space:
a direct application of Zorn Lemma to the family of all closed
non-empty $G$-subspaces of $X$ leads to the 
existence of a minimal $G$-subspace, $Y$,
and of course the restriction of the action of $G$ to $Y$ is still free.
Hence the following corollary, which was, in fact, the 
{\it raison d'\^etre} of Veech's result:

\begin{corol}
Every locally compact group admits an effective minimal action on a
compact space. \qed
\end{corol}

In particular, every locally compact group admits a
fixed point-free action on a compactum.

\begin{remark} As pointed out to me by A. Kechris, yet another proof
of the Veech Theorem can be found in \cite{AS}.
\end{remark}

\begin{remark}
It is rather surprising that beyond the class of locally compact
groups and some of their most immediate derivatives (such as MAP
groups, including free topological groups and additive
topological groups of locally convex
spaces) we know precious little
about topological groups admitting free actions on compacta.
One exception is the rather exotic class of {\it P-groups}
\cite{P21}, that is, those
topological groups in which every $G_\delta$-subset is open --- hardly
a class of great significance!
Are there any other visible classes of
topological groups admitting free actions on compacta? What we
rather have at
the moment, is a pageant of counter-examples growing richer by the day,
cf. below.
\end{remark}

\subsection{Fixed point on compacta property}
If one wishes to go all the way
in the opposite direction from the existence of
free actions, here is the concept to suit. One says that a
topological group $G$ has the {\it fixed point on compacta property}
({\it f.p.c.}),
or else is {\it extremely amenable,} if $G$ has a fixed point in
every compactum $X$ it acts upon: $g\cdot x^\ast=x^\ast$ for some
$x^\ast\in X$ and all $g\in G$.

According to the Veech Theorem, such topological groups can never be
locally compact (in particular, discrete). 
While it is not difficult to contruct even discrete
semigroups with the fixed point on compacta property 
\cite{Gra}, it was at first
unclear if topological groups with this property existed at all, as
is documented by the relevant question asked in print by 
T. Mitchell in 1970 \cite{Mit}. The examples of this kind
were hard to come by, and to the best of our knowledge, it was first done
by Herer and Christensen \cite{HC}
(even though they appeared to be unaware of Mitchell's question).

Before discussing such examples, let us make a few
fleeting remarks on amenable groups. 
A topological group $G$ is said to be
{\it amenable} if it possesses an {\it invariant mean}, that is,
a linear real-valued 
functional $\phi$ on the Banach space $C^b_\Rsh(G)$, 
which is positive (that is, $\phi(f)\geq 0$ whenever
$f$ is a non-negative function), of norm one, 
and invariant under the left action of
$G$, that is,
\[\phi(f) = \phi(g\cdot f) \mbox{ for all $g\in G$ and
$f\in C^b_\Rsh(G)$.}\]

The two most immediate classes of amenable groups are given by 
compact groups, where the invariant mean is just the Haar integral,
\[\phi(f):=\int_G f(x)d\mu(x),\]
and abelian topological groups, for which the proof is somewhat more
involved. 
The problematics of amenability has grown out of the famous
Banach--Tarski paradox (which essentially amounts to the 
non-amenability of the free groups on two generators). 
A good introduction to amenable groups is the book \cite{Wagon},
while the book \cite{Gre} is a classical reference, and the monograph
\cite{Pat} is the most modern and comprehensive source
with greater emphasis put on links with modern analysis.

Notice that every topological group with the fixed point on compacta
property is amenable. Indeed, every function $f\in C^b_\Rsh(G)$
extends to a unique continuous function $\bar f$
on the greatest ambit, and by evaluating it at a
chosen fixed point $x^\ast\in \gamma_\Rsh(G)$
one obtains an invariant mean,
\[\phi(f):=\bar f(x^\ast),\]
which is even multiplicative:
\[\phi(fg)=\phi(f)\phi(g) \mbox{ for all $f,g$.}\]
This explains the origin of the name `extremely amenable group.'

On the other hand, if $G$ is an amenable topological group admitting
no nontrivial unitary representations, then $G$ has the fixed point on
compacta property, that is, $G$ is extremely amenable. 
Indeed, suppose that $G$ does not have the f.p.c., that
is, admits a nontrivial minimal action on a compact space $X$. 
Fix a point $x_0\in X$, and let $h\colon G\mapsto g\cdot x_0$ be
the corresponding orbit map, $h\colon G\to X$.
The space $C(X)$ can be made into a pre-Hilbert space through
endowing it with the
(positive semi-definite) inner product:
\[(f,g):=\phi((h\circ f)(\overline{h\circ g})),\]
where $\phi$ denotes an invariant mean for $G$. The group $G$ acts
strongly continuously by isometries on the space $C(X)$, and this
representation factors through to the associated Hilbert space $\H$,
giving rise to a unitary representation. To prove that the
obtained unitary representation is
non-trivial, one uses the minimality of $X$.
(It is easy to show the existence of a non-zero function $f$ on $X$ 
whose support is disjoint from the support of a suitable
translation of $f$; then the $G$-orbit of the image of
$f$ in $\H$ is necessarily non-trivial.)
Thus we arrive at the historically first example of a topological
group satisfying f.p.c. property.

\begin{example}
It follows that the above mentioned example \ref{hererex}
by Herer and Christensen \cite{HC} of
an abelian topological group without nontrivial unitary representations
has the f.p.c. property. 
\end{example}

\begin{example}
The examples \ref{banaex} by
Banaszczyk also have the f.p.c. property, and thus the Veech
Theorem cannot be extended from locally compact groups even to 
Banach--Lie groups.
\end{example}

In all the fairness, the above two examples look more 
like genuine, elaborately designed {\it counter-}examples,
and even the name under which topological groups without unitary
representations appear in the above quoted papers 
-- {\it exotic} topological groups --
bears a witness to that.
However, more recent developments have revealed a highly surprising
trend: among `massive' groups, the fixed point on compacta property
is rather common! Here are some further examples known to date.

\begin{example} \fbox{$O(\H)_s$} Denote by $O(\H)_s$
the group formed by all orthogonal operators on the 
infinite dimensional separable real
Hilbert space $\H\cong l_2$, equipped with the strong operator
topology. 
Gromov and Milman have shown in 1983 \cite{GrM} that
${O(\infty)}$ has the fixed point on
compacta property. Similarly, the full unitary group $U(\H)_s$ of
the complex Hilbert space with the strong operator topology
has the fixed point on compacta property as well.
\label{gromov1}
\end{example}

\begin{example} \fbox{$L_1(X,{\mathbb T})$}
Let $(X,\mu)$ denote a non-atomic Lebesgue probability space.
(It is known that every two such spaces are isomorphic, 
and a standard model of a non-atomic Lebesgue space is
the closed unit interval $\mathbb I$ equipped with the usual
Lebesgue measure; a measure-preserving isomorphism between
subsets of full measure of $\mathbb I$ and $X$ is called a
{\it parametrization} of $X$.) Denote by $L_1(X,{\mathbb T})$ the
group formed by all measurable maps from $X$ to the circle rotation
group ${\mathbb T}= U(1)$, where the group operations are defined
pointwise. Make $L_1(X,{\mathbb T})$ into a
topological group using the $L_1$-distance between
measurable functions $f,g\colon X\to {\mathbb T}$:
\begin{equation}
d_1(f,g):=\int_X d(f(x),g(x))d\mu(x).
\end{equation}
Here $d$ is any metric on the circle group.
Eli Glasner \cite{Gl} (and, independently, Furstenberg and B. Weiss,
unpublished) have proved that the topological group 
$L_1(X,{\mathbb T})$ has the fixed point on compacta property.
\label{glasner1}
\end{example}

\begin{example}
\fbox{$\Homeo_+({\mathbb I})$ and $\Homeo_+({\mathbb R})$}
\label{pestov}
The present author has proved \cite{P19} that the groups
of orientation-preserving homeomorphisms of the closed interval,
$\Homeo_+({\mathbb I})$, and of the real line, $\Homeo_+({\mathbb R})$,
equipped with the compact-open topology, have the fixed point on
compacta property.
\end{example}

\begin{example} \fbox{$\Aut(X,\mu)$}
\label{giordano}
Again, let $X=(X,\mu)$ be a non-atomic Lebesgue measure space.
An {\it automorphism} (or a {\it measure-preserving transformation})
of such a space is a measurable invertible map $f$ from $X$ to itself
such that for every measurable $A\subseteq X$ one has
$\mu(f^{-1}(A))=\mu(A)$. (Then the inverse map $f^{-1}$ is
automatically measurable as well.) The collection of all automorphisms
of $X$ forms a group, $\Aut(X)$. The {\it strong} topology on the
group $\Aut(X)_u$ is the Hausdorff group topology whose
neighbourhood basis at the identity is formed by all sets of the form
\[N(Y,\e):=\{h\in\Aut(X)\colon \mu(Y\cap h(Y))<\e\},\]
where $\e>0$ and $Y\subseteq X$ is measurable. The {\it uniform}
topology on $\Aut(X)$ is a (finer) group 
topology generated by the bi-invariant metric 
\begin{equation}
d_{unif}(g,h):=\mu\{x\in X\colon g(x)\neq h(x)\}.
\label{uniform}
\end{equation}
(Cf. e.g. \cite{Hal}, pp. 65 and 69--74).

The recent joint result by Thierry Giordano and the present
author \cite{GP} says that the group $\Aut(X)$ with the strong topology
has the fixed point on compacta property.
(This result admits further generalizations and ramifications,
which will be hopefully explored in one or more papers 
by the same authors, currently in preparation.)
\end{example}

In our view, the entire trend is very significant, indicating that
the properties of massive groups in some respect are completely opposite
to those of locally compact groups. Moreover, the way the majority 
of the above
results are being established provides an opportunity to link
theory of topological groups with an important
development in modern analysis and geometry --- the {\it phenomenon of
concentration of measure on high-dimensional structures.}

\subsection{Concentration of measure}
Let us describe the basic idea of concentration phenomenon, which is
meaningful in the presence of some sort of
proximity between points (usually distance, or
else uniformity) and `size' of sets (measure). We will adopt the
setting for analyzing concentration proposed by Gromov and
Milman in 1983 \cite{GrM}, which is provided
by a metric space $\Omega=(\Omega,\rho)$
equipped with a normalized ($\mu(\Omega)=1$) positive Borel measure.
(It should be remarked that this setting is, in all the likelihood, not
final, and at least two alternative frameworks for the concentration
phenomenon have been proposed recently:
an `ergodic' setting by Gromov \cite{Grogafa}
and an `affine' setting, cf. the preprint 
by Giannopoulos and Milman \cite{GiM}.)

Suppose $\Omega=(\Omega,\rho,\mu)$ is a metric space equipped with  
measure as above (mm-{\it space}).
Let $A\subseteq \Omega$ be an arbitrary Borel subset containing at least half
of all points, that is, $\mu(A)\geq \frac 12$. Denote by
${\mathcal O}_\e(A)$ the open $\e$-neighbourhood of $A$ in $\Omega$. 
How massive is ${\mathcal O}_\e(A)$? Or, equivalently, how small
--- in the sense of measure --- is
the `cap' $X\setminus {\mathcal O}_\e(A)$? (Fig. 3.)

\bigskip
\begin{center}
\scalebox{0.3}[0.3]{\includegraphics{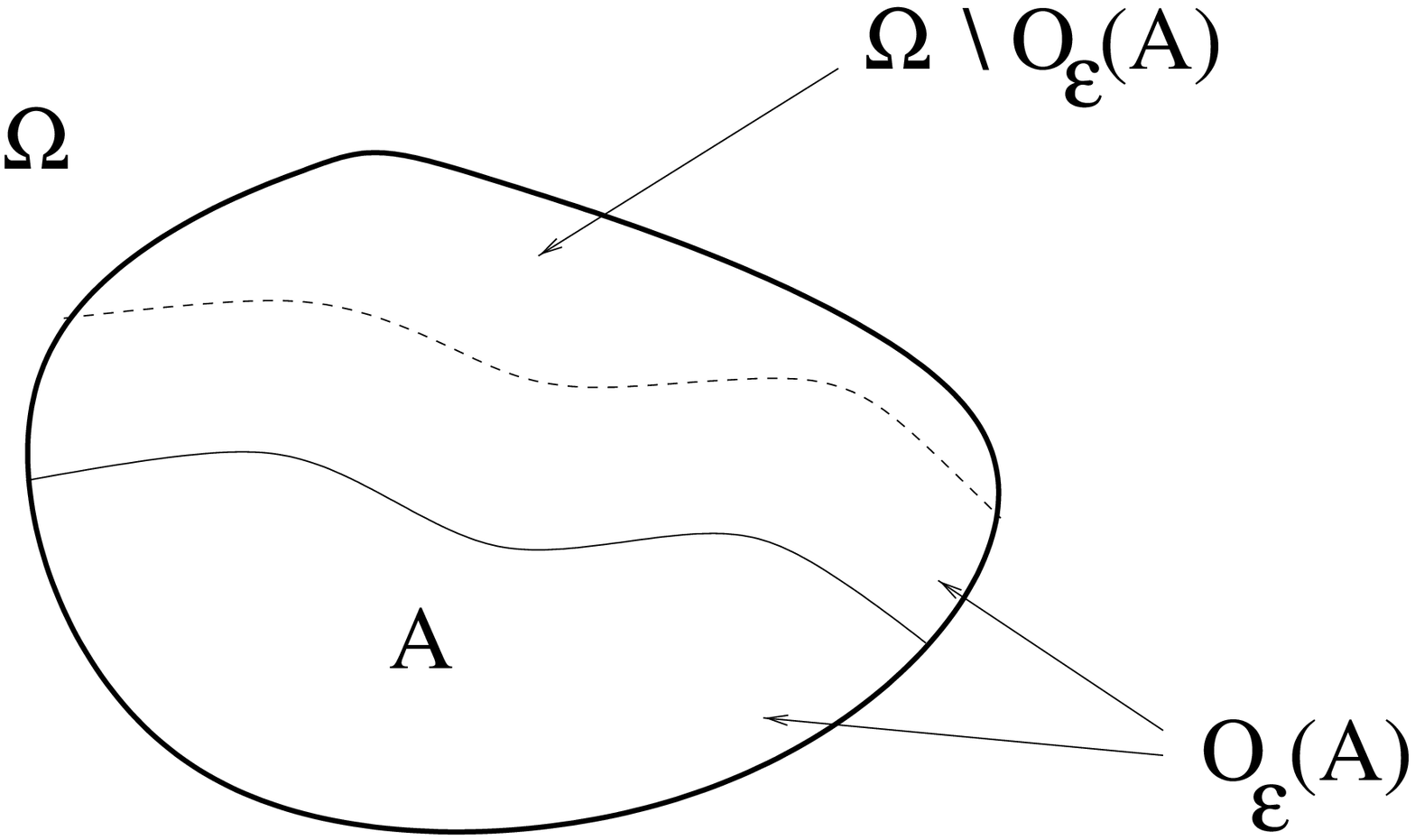}}
\smallskip

{\small
Fig. 3. An illustration to the concentration of measure.}
\end{center}

\bigskip

If $\Omega={\mathbb I}$ is the closed unit interval
equipped with the usual distance and the Lebesgue measure, then, by
letting $A=[0,\frac 12]$, one can see that the cap (which is, in this
case, the interval $(\frac 12+\e,1]$) need not be really
small in that it has
measure $\frac 12-\e$. Things do change however when we proceed to
higher-dimensional objects. Here is a heuristic way to describe 
the phenomenon of concentration of measure on structures of high
dimension: 
\smallskip

{\it if $\Omega$ is `high-dimensional' then, typically, 
the size of the `cap'  $\Omega\setminus {\mathcal O}_\e(A)$ is extremely
close to zero already for small values of $\e>0$. }
\smallskip

In other words,
nearly all points of $\Omega$ are $\e$-close to $A$ provided
$\mu(A)\geq\frac 12$.

A convenient way to quantify the concentration phenomenon
is to consider the {\it concentration function} of $\Omega$ which gives
the least upper bound on the measures of all `caps' as above: 
\begin{equation} 
\alpha_\Omega(\e)=1-\inf\left\{\mu\left({\mathcal O}_\e(A)\right)
\colon A\subseteq \Omega \mbox{ is Borel and }
\mu(A)\geq \frac 1 2\right\}.
\end{equation}
For the unit interval $\alpha(\e)=\frac 12-\e$, which is not very
interesting. However, things do look different if we consider, for
example, the $n$-spheres
$\s^n=\{x\in\R^{n+1}\colon \norm x_2 =1\}$, equipped with the 
geodesic distance and the 
(unique) normalized ($\mu_n(\s^n)=1$) rotation-invariant measure,
$\mu_n$
(which for $n=1$ turns into the
arc length, for $n=2$ into the surface area, for
$n=3$ the volume element, and so forth.) 
The maximal size of the `cap' is achieved for 
$A=\s^n_-$, the hemisphere. (This is one of the equivalent forms
of the so-called {\it isoperimetric inequality.}) 
Now pretty straightforward calculations at the level of a good
first-year
calculus student enable one to compute the concentration functions
$\alpha_{\s^n}$
of the $n$-spheres. Here are their graphs in
dimensions $n=3,10,100,500$. (Fig. 4.)
\bigskip

\begin{center}
\scalebox{0.6}[0.6]{\includegraphics{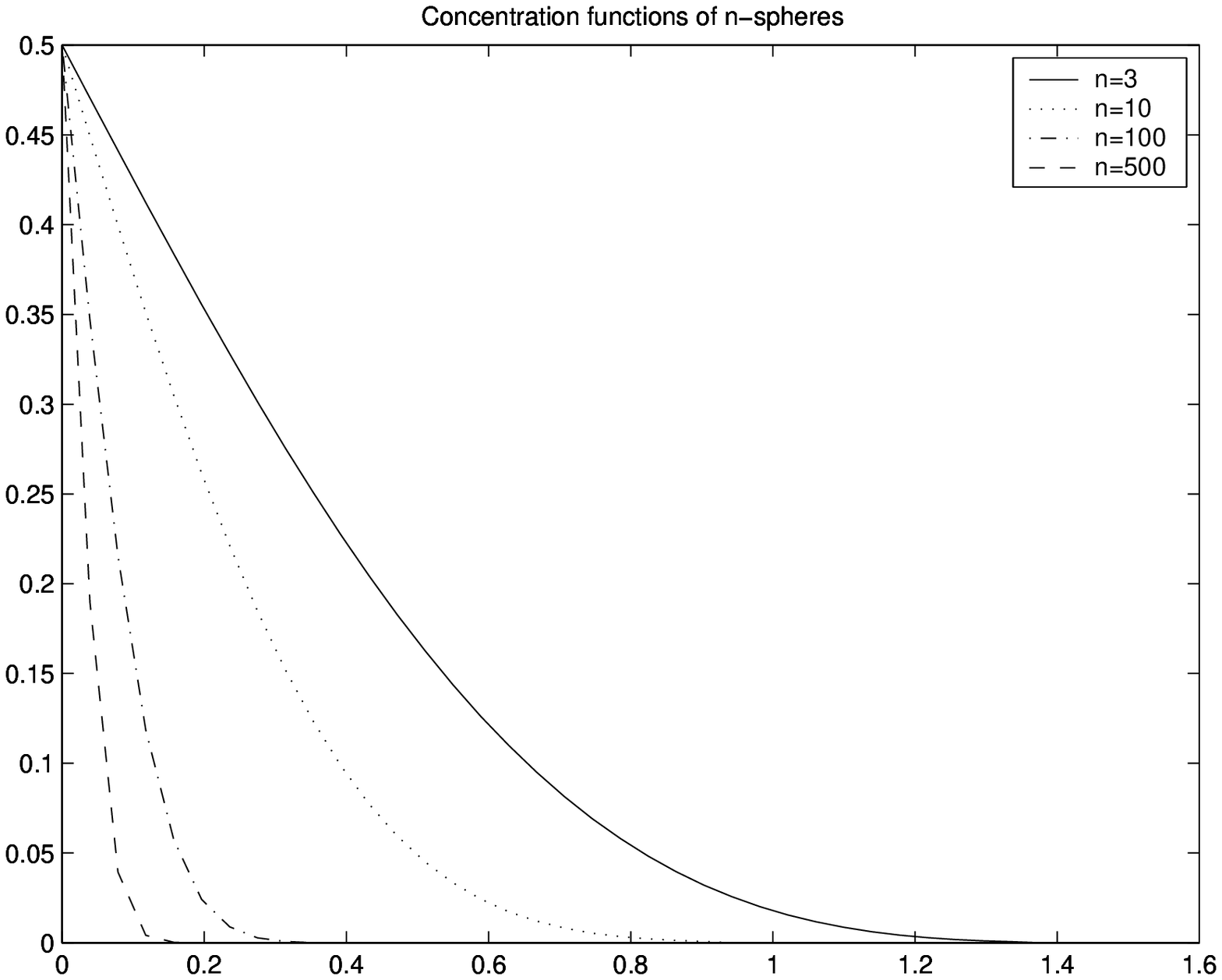}}
\smallskip

{\small
Fig. 4. Concentration functions of $n$-spheres,
$n=3,10,100,500$.}
\end{center}
\bigskip

The following result is simple but important. In particular, it
explains the origin of the terminology (`concentration of measure.')

\begin{thm}
Let $\Omega=(\Omega,\rho,\mu)$ be a metric space equipped with
a normalized Borel measure, and let $f\colon \Omega\to \R$ be 
a uniformly continuous function. Denote by $\delta=\delta(\e)$ the
modulus of uniform continuity of $f$, that is, for every
$x,y\in \Omega$ and $\e>0$ one has $\vert f(x)-f(y)\vert<\e$ whenever
$\rho(x,y)<\delta(\e)$. Denote by $M$ the median value of $f$
on $\Omega$, that is, a real number with 
\[\mu\{x\in\Omega\colon f(x)\leq M\} =
\mu\{x\in\Omega\colon f(x)\geq M\}.\]
Then the set of all $x\in\Omega$ such that
$\vert f(x)-M\vert<\e$ has measure at least
\[1 - 2\alpha_\Omega(\delta(\e)).\]
\qed
\label{median}
\end{thm}

If the concentration function $\alpha$ drops off sharply for small
values of the argument, then most of the points of the domain
$\Omega$ `concentrate' near one value of $f$. In other words,
the function $f$ is, from the probabilistic viewpoint, almost constant.

The asymptotic behaviour of families of spaces with metric and
measure, namely the tendency of concentration functions to
fall off sharply near zero as the dimension grows (evident in Fig. 4), 
can be formalised as follows.
One says that an infinite family
$(\Omega_n,\rho_n,\mu_n)$ of metric spaces equipped with measure
is a 
 {\it normal L\'evy family} if
\[\alpha_{\Omega_n}(\e)\leq C_1 e^{-C_2\e^2n}.\]
In simpler words, it means that the measures of `caps' go to zero
exponentially fast in dimension for a fixed value of $\e>0$.
For example, $n$-spheres form a normal L\'evy family. 
What is more, and this is very important,
`naturally occurring' infinite families of probabilistic metric
spaces are, typically, normal L\'evy.
Here are just three examples important for what follows.

\begin{example} (Glasner \cite{Gl}, and, independently, 
Furstenberg and B. Weiss.)
The family of tori ${\mathbb T}^n$, equipped with
the normalized Haar measure and the metric
\[d(x,y):=\frac 1 n\sum_{i=1}^n \vert x-y\vert,\]
where the absolute value is induced through the standard
embedding ${\mathbb T}\equiv U(1)
\subset\C$, form a normal L\'evy family.
(This follows from more general results of Talagrand \cite{Ta}.)
\label{glasner2}
\end{example}

\begin{example}
(Maurey, \cite{Ma}.)
The groups of permutations $S_n$ of rank $n$, equipped
with the normalized Hamming distance
\[d(\sigma,\tau):=\frac 1n \abs{\{i\colon \sigma(i)\neq\tau(i)}\]
and the normalized counting measure
\[\mu(A):=\frac{\abs{A}}{n!}\]
forms a normal L\'evy family, with the concentration functions
satisfying the estimate
\[\alpha_{S_n}(\e)\leq \exp(-\e^2n/64).\]
\label{maurey2}
\end{example}

\begin{example}
(Gromov and Milman, \cite{GrM}.) The special orthogonal groups 
$SO(n)$ of rank $n$ consist of all orthogonal $n\times n$ matrices with
real entries having determinant $+1$. The family of
these groups, equipped with the normalized 
Haar measure and the uniform metric
(that is, the metric induced by the operator norm under the
standard embedding $SO(n)\subseteq {\mathcal L}(\R^n)$), form a
normal L\'evy family.
\label{gromov2}
\end{example}

Listed below are just a few common manifestations of the phenomenon of
concentration of measure in mathematical sciences.\smallskip

$\bullet$ The Law of Large Numbers: the average value of a long
sequence of $0$s and $1$s obtained by tossing a fair coin
is typically $\approx\frac 12$.

$\bullet$ Most of the volume of a high-dimensional 
Euclidean ball is concentrated
near the surface.

$\bullet$ Most of the volume of a high-dimensional unit cube is concentrated
near the corners.

$\bullet$ Blowing-Up Lemma in coding theory: if the Hamming cube
$\{0,1\}^n$ is partitioned into two subsets of equal size, then almost all
binary $n$-strings are close to both subsets.

$\bullet$ Dvoretzky Theorem:
If a convex body in a high-dimensional space is cut by a random
plane, the section typically looks almost like a circle.

$\bullet$ Two random vectors in a high-dimensional Euclidean space
selected independently of each other are typically nearly orthogonal. 
\smallskip

Various aspects of the phenomenon of concentration of measure on
high-dimensional structures, including all the above
examples, are discussed in 
\cite{GrM}, \cite{M1}, \cite{M2}, \cite{MS}, \cite{Ta},
\cite{GiM}, \cite{Grogafa}, and \cite{Gro}, Ch. 3$\frac 12_+$. 

\subsection{L\'evy groups}
In topological algebra the concentration phenomenon is
captured by the concept of a L\'evy group. The definition below
slightly extends the original one \cite{GrM} (cf. \cite{Gl}
and \cite{P19}) in that the metric is replaced with uniformity.

\begin{defin}
We say that
a topological group $G$ is a {\it L\'evy group} if there is a 
family $\mathcal K$ of compact subgroups of $G$ with the following
properties.
\begin{enumerate}
\item The family $\mathcal K$ is directed by inclusion, that is,
for any $F,H\in{\mathcal K}$ there is a $K\in {\mathcal K}$
with $F\cup H\subseteq K$. 
\item The union $\cup{\mathcal K}$ is everywhere dense in $G$.
\item Let a family of Borel 
subsets $A_K\subseteq K$, $K\in{\mathcal K}$
have the property that 
\begin{equation}
\liminf_{K\in{\mathcal K}}\mu_K(A_K)>0,
\label{assumption}
\end{equation}
where $\mu_K$ denotes the normalized Haar measure on $K$. Then for
every neighbourhood of zero, $V$, in $G$,
\begin{equation}
\lim_{K\in{\mathcal K}}\mu_K(K\cap (VA_K))=1.
\label{conclusion}
\end{equation}
\end{enumerate}
\end{defin}

\begin{remark} Zero on the r.h.s. of 
(\ref{assumption}) can be replaced, without any loss in generality, 
by any positive constant $<1$, for example $\frac 12$. 
\end{remark}

\begin{examples} 1. The group $L(X,{\mathbb T})$ (Ex. \ref{glasner1})
forms a L\'evy group. 
Simple functions, constant on
elements of a sequence of refining partitions of the Lebesgue space
$X$, form an increasing sequence of tori having everywhere dense
union in the group. Now one applies the observation from Ex. \ref{glasner2}.
\smallskip

2. The group ${O(\H)}_s$ with the strong operator topology 
(Ex. \ref{gromov1}) is a L\'evy group. It follows from 
Ex. \ref{gromov2} and the following observation.
Let us identify elements of $SO(n)$ with those orthogonal
operators in $\H$ which are
represented, with respect to a chosen orthonormal basis in $\H$,
by matrices with only finitely many non-zero entries. 
Then the union of 
the increasing sequence of the special
orthogonal groups of growing finite rank embedded into
each other via
\begin{equation}
SO(n)\ni A\mapsto \left(\begin{matrix}
1 & 0_{1\times (n-1)} \\ 0_{(n-1)\times 1} & A
\end{matrix} \right)  \in SO(n+1)
\end{equation}
is everywhere dense in ${O(\H)}_s$.

A similar argument applies in the case of infinite unitary groups.
\smallskip

3. The group $\Aut(X)$ of measure-preserving transformations of a
Lebesgue space equipped with the strong topology
(Ex. \ref{giordano}) is L\'evy. Indeed, consider a parametrization of
$X$ by the closed unit interval with the Lebesgue measure.
Call a transformation of the interval $\I$
a {\it permutation of rank
$n$} if it maps each binary interval of rank $n$ to such an interval
by a translation. Then it is well-known in ergodic
theory (a corollary of Rokhlin's Lemma) that the collection of all
such permutations is everywhere dense in $\Aut(X)$ (the so-called
Weak Approximation Theorem, cf.  e.g. \cite{Hal},
pp. 65--68). The uniform metric on $\Aut(X)$ induces the normalised
Hamming distance on each group of permutations $S_n$, and by using
Ex. \ref{maurey2}, one concludes that $\Aut(X)$ is a L\'evy group.
\end{examples}

In order to establish the fixed point on compacta property
for the groups from Examples \ref{gromov1}, \ref{glasner1},
\ref{giordano}, it is now sufficient to establish the following.

\begin{thm} Every L\'evy group has the fixed point on compacta property.
\label{gm}
\end{thm}

This result belongs to Gromov and Milman (\cite{GrM}, Th. 5.3), 
who stated it in a somewhat more restricted form, 
later removed by Glasner \cite{Gl}, Th. 1.2. 
See also \cite{P19}, Th. 9.1.

Let us prove Theorem \ref{gm}.
It is clearly sufficient to establish the existence of a fixed
point for the canonical action of $G$ on the greatest 
ambit $\gamma_\Rsh(G)$. Furthermore, the
compactness considerations easily imply that it is enough
to find the common fixed point for
an arbitrary finite subset of elements $g_1,g_2,\dots,g_n$
of $G$. Such a fixed point
will certainly exist if the following property is
satisfied: for every element $V$ of the unique uniform structure
on $\gamma_\Rsh(G)$, there is a point $x\in \gamma_\Rsh(G)$ 
such that all he elements
$g_1,g_2,\dots,g_n$ fail to move $x$ beyond the neighbourhood
$V[x]$. One does not lose in generality by assuming that
$x\in G$, and instead of the entourage $V$ one can consider,
using a common trick in uniform topology, an 
arbitrary bounded ${\mathcal U}_\Rsh(G)$-uniformly continuous
function $f$ from $G$ to a finite-dimensional Euclidean space.
The property we want to establish 
becomes this: for every such $f$ as above and every $\e>0$
there is an $x\in G$ such that for all $i$, 
\[\abs{f(g_ix)-f(x)}<\e.\]

Assume for the reasons of mere technical simplicity that
the L\'evy group $G={\cup_{i=1}^\infty}G_i$ under
consideration is separable,
so that a net of approximating subgroups can be replaced with
an increasing chain, whose union in addition coincides with $G$, 
and also that $G$ is metrizable, and fix a
right-invariant metric $\rho$ generating the topology of $G$.
\par
We denote by $\mu_i$ the normalised Haar measure on the compact
group $G_i$.
The elements $g_1,g_2,\dots,g_n$ from the given 
finite collection are contained in a $G_N$ for $N$ sufficiently large,
and thus we can assume by removing the first $N-1$ groups in the
sequence that $g_1,g_2,\dots,g_n\in G_1$.
Let $f$ be a ${\mathcal U}_\Rsh(G)$-uniformly continuous
bounded function on $G$ taking values in a finite-dimensional
Euclidean space $\R^k$. We will consider the $l_\infty$-norm on the 
latter space, just to choose any.
Denote by $\delta=\delta(\e)$ the modulus
of continuity of $f$. Let $f_j$, $j=1,\dots,k$ be the
components of $f$, and denote for each $i\in\N$ by $M_{i,k}$
the median value of the restriction $f_k\vert_{G_i}$. 
Let $\e>0$ be any.
According to Theorem \ref{median}, for all 
elements $g\in G_i$
from a set of measure at least $1-2\alpha_{G_i}(\delta(\e))$
one has 
\[\abs{f_{i,k}(g)-M_{i,k}}<\e.\]
Replacing $\e$ with $2\e$,
using the compactness of the closed interval, and proceeding to a
subsequence if necessary, 
one can assume that the numbers
$M_{i,k}=M_k$ are independent of $i$. Let $M=(M_1,M_2,\dots,M_k)\in\R^k$.
It follows that for each $i\in\N$, for all $g\in G_i$ from a set $S_i$
of measure at least $1-2\alpha_{G_i}(\delta(\e))$ one has
\[\norm{f(g)-M}_\infty<\e.\]
The translates of the set $S_i$ by elements $g_1^{-1},\dots,g_n^{-1}$
have the same measure as $S_i$ and the intersection
of all such translates, which we will denote $X_i$, is of measure 
(computed in the group $G_i$) at least
$1-2n\alpha_{G_i}(\delta(\e))$. The definition of a L\'evy group
means that the concentration functions
$\alpha_{G_i}$ converge pointwise to zero as $i\to\infty$. In particular,
for $i$ sufficiently large, the sets $X_i$ are of positive measure.
Any point $x\in X_i$ will then have the desired property. \qed

\begin{remark}
The proof of the f.p.c. property of the groups from example
\ref{pestov} (that is, the homeomorphism groups of the closed
and the open unit intervals) is established
in \cite{P21} in a somewhat different
manner, using the infinite Ramsey theory, and we do not reproduce
it here, cf. e.g. our short survey
\cite{P19}. Interestingly,
it has been repeatedly noted (cf. e.g. \cite{M2})
that the Ramsey theorems in combinatorics are very close
in spirit to the phenomenon of concentration of measure. 
\end{remark}

The size of the class of extremely amenable topological
groups is immense. The following result
(conjectured independently by Gromov and Uspenskij in private
discussions with the author) 
was recently established by the present author.

\begin{thm}[Pestov \cite{Pnew}]
Let $U$ be an $\omega$-homogeneous generalized Urysohn metric space.
Then the group $\Iso(U)$ has the fixed point on compacta property.
\label{urr}
\qed
\end{thm}

The proof explores both the phenomenon of 
concentration of measure and a close link
between the extreme amenability of the groups of isometries
$\Iso(X)$ of sufficiently homogeneous metric spaces $X$ and a
Ramsey-type property of $X$.

Theorem \ref{urr} and Uspenskij's Theorem \ref{uspee}
together imply:

\begin{corol}
Every topological group embeds, as a topological subgroup, into
an extremely amenable topological group, that is, a topological
group with the fixed point on compacta property. \qed
\label{toanex}
\end{corol}

Even if one replaces `extremely amenable' with 
`amenable,' the result remains new. Notice
that amenability is inherited by topological subgroups of
{\it locally compact} amenable groups; for non-locally compact
groups this is no longer true \cite{dlH}, and Corollary 
\ref{toanex} takes this observation to its extreme.

Finally, here is another corollary of Theorem \ref{urr}, 
answering a question from \cite{U1}.

\begin{corol}
The topological groups
$\Iso(\U)$ and $\Homeo(\I^\omega)$ are not isomorphic.
\end{corol}

\begin{proof}
Indeed, the latter group admits a continuous action without fixed
points on the compact space $\I^\omega$.
\end{proof}

Thus, the two examples of universal Polish groups we are aware of
are different indeed.

\subsubsection*{Some further questions on 
concentration in topological groups}
The examples of L\'evy groups belonging to Gromov and Milman
($U(\H)_s$, Ex. \ref{gromov1}) and to Glasner and Furstenberg--B. Weiss
($L(X,U(1))$, Ex. \ref{glasner1}) 
share in fact a profound similarity in that
both of them are unitary groups of suitable von Neumann algebras
equipped with the ultraweak topology
\cite{Sa}: $L_\infty(X)$ in the
first case, $L({\mathcal H})$ in the second.

Moreover, if a von Neumann
algebra $W$ is such that the unitary group with the ultraweak topology
is L\'evy (and therefore, in particular, 
amenable), then $W$ is hyperfinite (de la Harpe \cite{dlH} and Paterson
\cite{Pat}) and therefore injective.
The following question is natural.
 
\begin{problem} What are those von Neumann algebras
whose unitary groups with the ultraweak topology are L\'evy?
\end{problem}

\begin{problem}
The author understands that (a version of)
the following problem
was put forward by Furstenberg at least 17 years ago:
{\it Let $G$ be a topological group that is the union
of a directed family of compact subgroups. Is it possible to
express the property of $G$ being L\'evy through the
existence of fixed points in some compacta that $G$ acts upon?}

If one interprets the problem as whether or not the L\'evy
property for topological groups of the above type is equivalent to
extreme amenability, then we strongly suspect that the answer is no,
though at the moment we do not have any concrete
counter-example.
However, it is conceivable that the problem admits a wider
interpretation, leading to a positive answer.
\end{problem}

\subsection{Extreme amenability and left syndetic sets}
It is worth stressing that though the fixed point on compacta
property is formulated in exterior terms (actions on compact spaces),
it is an intrinsic property of a topological group $G$ itself.
This becomes evident if one looks at the following alternative
criterion, which can also be used to establish the fixed point theorems.
Recall that a subset $S$ of a group $G$
is {\it left syndetic,}
or ({\it left}) {\it relatively dense,} 
if $FS=G$ for some finite subset $F\subseteq G$.

\begin{thm}[Pestov, \cite{P19}] A topological group 
$G$ has the f.p.c. property if and only if
for each left syndetic 
subset $S\subset G$, the set $SS^{-1}$ is everywhere dense
in $G$.
\label{criterion}
\qed
\end{thm}

Theorem \ref{criterion} was inspired by, and is to be compared
with, the following classical result. 

\begin{thm}[F\o lner \cite{Fol}; Cotlar---Ricabarra \cite{CotRic}; 
Ellis---Keynes \cite{EK}]
An abelian topological group is minimally almost periodic 
if and only if for each big $S\subseteq G$, the set
$S-S+S$ is everywhere dense in $G$.
\label{fcrek}
\qed
\end{thm}

To better appreciate the similarity between \ref{criterion}
and \ref{fcrek}, notice
that every abelian topological group with the fixed point on compacta
property is minimally almost periodic. 

It is in fact unknown if the converse is true!
Since the f.p.c. property intuitively feels 
so much stronger a restriction than minimal almost
periodicity, our inability to distinguish between the two
properties comes as a surprise.

Both results can be turned the other way round. In particular,
the following mirror image of
Theorem \ref{fcrek} yields a criterion
for the existence of sufficiently many characters. 

\begin{corol}
An abelian topological group $G$ is maximally almost periodic if and
only if for every $g\in G$, $g\neq 0$, there exists a big set
$S\subseteq G$ such that the closure of $S-S+S$ does not contain $g$.
\end{corol}

In fact, this is a sort of result that gives a fair idea of what 
would be an
acceptable answer to Shtern's question \ref{questshtern}.

It is natural to ask for a non-abelian version of the above, with
$S-S+S$ being replaced by $S^{-1}S^2S^{-1}$.
While the answer seems to be unknown in the full generality, an important
advance is due to Landstadt \cite{Land}
who established the result for amenable
topological groups. 

The following particular case of Theorem \ref{fcrek} is of
a special interest in combinatorial number theory.

\begin{corol}
If $S$ is a relatively dense subset of the integers, then
$S-S+S$ is a neighbourhood of zero in the Bohr topology on the group
$\Z$.
\end{corol}

It remains unknown for a long time \cite{V2} if one can replace
in the above result $S-S+S$ with $S-S$. 

Glasner \cite{Gl} 
has observed that a negative answer would follow if one
constructs an example of a minimally almost periodic, monothetic 
topological group without the f.p.c. property. For a simpler 
explanation of why, see also \cite{P19}.
No such example is presently known.
To construct it, one apparently needs to maintain a
very fine balance between minimal
almost periodicity and a property that goes in 
exactly the opposite direction to measure concentration: it is 
some form of {\it measure dissipation,} cf. \cite{Gro}.

\section{Parallels between topological and discrete groups}
\subsection{Subgroups of finitely generated groups}
\subsubsection{Higman--Neumann--Neumann theorem}
Here is an example of how actions can be used as a tool 
in theory of topological groups.
Consider the following classical result in group theory
(no topology present!).

\begin{thm}[Higman--Neumann--Neumann, \cite{HNN, NN}]
Every countable group is isomorphic with a subgroup of a
2-generated group. \qed
\label{hnn}
\end{thm}

Somewhat unexpectedly, the above result has a direct
counterpart for topological groups.

\begin{thm}[Morris and Pestov, \cite{MP1}]
\label{mp}
Every countable topological group is isomorphic with a
subgroup of a group algebraically generated by two elements.
\end{thm}

The following corollary is more or less straightforward and
puts the result in a natural topological group wrapping.

\begin{corol}[\cite{MP1}]
Every separable topological group is isomorphic with a
topological subgroup of a group with two topological generators.
\qed
\label{sep}
\end{corol}

From here one can deduce without much effort a description of
topological subgroups of topologically finitely generated 
groups. A topological group is
called $\omega$-{\it bounded} if it is covered with countably
many translations of every non-empty open subset. This concept is
modelled on that of a totally bounded group. A group is
$\omega$-bounded if and only if it embeds into the direct product of
a family of separable metrizable groups (with the usual product
topology).

\begin{corol}[\cite{MP1}]
A topological group embeds into a topologically
finitely generated topological group
if and only if it is $\omega$-bounded and has weight at most
continuum. \qed
\label{iff}
\end{corol}

The proof of Theorem \ref{mp} consists of nothing more than injecting
a bit of topological dynamics into a proof 
of the Higman--Neumann--Neumann Theorem \ref{hnn}
due to Galvin \cite{gal}.
Here is its outline. Enumerate the group in question with odd
positive integers,
$G=\{g_1,g_3,\dots, g_{2k+1},\dots\}$. 
Let $X$ be a set whose group of permutations contains
$G$: $G\hookrightarrow \Aut(X)$. Form a new set
\[\tilde X=\Bbb Z\times\Bbb Z\times X \cong
\oplus_{(m,n)\in\Bbb Z\times\Bbb Z}\{(m,n)\}\times X\}.\]
Now define permutations $a$ and $b$ of 
$\tilde X$ by letting (1)
$a\cdot (m,n,x)=(m+1,n,x)$, (2)
$b\cdot (0,n,x)=(0,n+1,x)$, (3) 
$b\cdot (m,n,x)= (m,n,xg_m)$
 if $m$ is odd, $m>0$, and $n\geq 0$, and finally (4)
making $b$ leave $(m,n,x)$ fixed otherwise.

Embed $G$ into the group of permutations of $\tilde X$ by
letting it act on $\{0,0\}\times X$ in a way identical to its action on $X$,
and on the rest of $\tilde X$ in a trivial way (every point is fixed).
Straightforward computations show that $G$ is contained in the group
${\mathrm{gp}}(a,b)$ generated by permutations $a$ and $b$.
The Higman--Neumann--Neumann Theorem \ref{hnn} is thus proved.

To obtain from the above a proof of Theorem \ref{mp}, it
suffices to put on $X$ a compact topology so as to make $G$ into a
topological subgroup of $\Homeo(X)$ (Teleman's theorem!), and
to topologize $\tilde X$ as the disjoint sum of compacta. Both $a$
and $b$ are now homeomorphisms, and the compact-open
topology makes $\Homeo\tilde X$ into a topological group obviously
containing $G$ as a topological subgroup. Q.E.D.

\subsubsection{Subgroups of monothetic groups}
What about 1-generated or, as they are most commonly called,
{\it monothetic} topological groups? They are certainly abelian, and
so are all their subgroups. Nevertheless, this turns out to be the
only additional restriction one has to impose on all potential subgroups
of such groups.

\begin{thm}[Morris and Pestov, \cite{MP2}]
Every separable abelian topological group is isomorphic with a
topological subgroup of a group with one topological generator.
\label{ab}
\end{thm}

One can abelianize Corollary \ref{iff} as well.

\begin{corol}[\cite{MP2}]
A topological group embeds into a monothetic topological group
if and only if it is abelian, $\omega$-bounded and has weight at most
continuum. \qed
\label{iffab}
\end{corol}

What is manifest, is that the proof cannot be aped
after a discrete case, simply because Theorem \ref{ab} clearly has
no discrete counterpart! Instead, the proof is based on an
entirely different technique, giving us an opportunity to introduce
into consideration free (abelian) topological groups.

\subsubsection{Free (abelian) topological groups and free locally
convex spaces}
Let $X$ be a completely regular $T_1$ topological space. A topological
group $F(X)$ is called the
{\it free topological group} on $X$ if it contains a
topological copy of $X$ as a
distinguished topological subspace in such a way that the following
diagram can be made commutative for every continuous mapping $f$ from
$X$ to an arbitrary topological group $G$ by means of a unique
continuous homomorphism $\bar f$:
\begin{eqnarray*}
X & \hookrightarrow & F(X) \\
& \forall f \searrow  & \downarrow {\exists \bar f}\\
& & G
\end{eqnarray*}
In a completely similar way, one defines the
free {\it abelian} topological group, $A(X)$, and the 
{\it free locally convex space}, $L(X)$. (In the latter case,
the morphisms are continuous linear operators between locally
convex spaces.) 

One of the most immediate observations concerning all three
types of objects is their algebraic freedom: $F(X)$ and $A(X)$
are algebraically the free and the free abelian groups on $X$
correspondingly, while $L(X)$ is a vector space spanned by
$X$ as an algebraic (Hamel) basis. \par
The topology on both $A(X)$ and $L(X)$ can be easily described
using the following construction going back to Graev \cite{Graev} and
Arens--Eells \cite{AE}, see also \cite{Rai,Fl1, Fl2}.
For a pseudometric $\rho$ on the set $X^\dag=X\cup\{0\}$ denote
by $\bar\rho$ the maximal translation invariant 
pseudometric on $A(X)$ with the property that
$\bar\rho\vert X^\dag=\rho$. The existence of such a pseudometric is
obvious, and moreover its value can be computed rather explicitely,
or at least in combinatorial terms of manageable complexity. 
Now it is a matter of an easy exercise, to show that the collection of
all pseudometrics of the form $\bar\rho$ (called {\it Graev,}
or {\it maximal, pseudometrics}) determines the topology of the
free abelian topological group $A(X)$ as $\rho$ runs through the
collection of all compatible pseudometrics on $X\cup\{0\}$.
In a similar vein, let $p_\rho$ denote the 
maximal seminorm on $L(X)$ such that $p_\rho(x-y)=\rho(x,y)$,
$x,y\in X^\dag$. Again, the existence of the seminorm $p_\rho$ is
rather straightforward, and the seminorms of this form determine
the topology of $L(X)$ if one lets $\rho$ run over all compatible
pseudometrics on $X\cup\{0\}$. (Warning: in both cases,
it is not enough for $\rho$ to go
through {\it some} collection of generating pseudometrics for the topology
of $X$!)

\begin{remark}
 For the (non-abelian) free topological group 
$F(X)$, no simple description of topology similar to the above
is known. The explicit constructions of generating
pseudometrics, such as those proposed first by Tkachenko \cite{Tk3}
and then others, are pretty hard to work with. However,
a similar description is known for the so-called 
{\it free SIN group,} $F_{\mathrm{SIN}}(X)$ (also known as the {\it free
balanced group.}) The definition of course mimicks that of the
free topological group, where it is assumed that all topological
groups under consideration are SIN (= have the left and the right
uniform structures coincide). The group 
$F_{\mathrm{SIN}}(X)$ is also algebraically free over $X$, and its topology
is described by the family of all {\it Graev pseudometrics} on the group
$F(X)$. A Graev pseudometric in this case is defined as the maximal
bi-invariant pseudometric $\bar\rho$ whose 
restriction to $X^\dag:= X\cup\{e\}$ coincides
with the given pseudometric $\rho$. Only in exceptional cases do
the topologies of $F(X)$ and of $F_{\mathrm{SIN}}(X)$ coincide.
\end{remark}

There is no single comprehensive reference 
to the up-to-date theory of free topological groups. Some pointers 
can be found in \cite{Arh1, Arh2, CHR, P12}.

Listed below (in no paticular order)
are results about free (abelian) topological groups 
that the present author likes most. (He hopes to be forgiven for
including among them $1\frac 13$ results of his own.)

\begin{itemize}
\item The semi-classical work by Markov \cite{Mar} and 
Graev \cite{Graev}.
\item Arhangel'ski\u\i's results on zero-dimensionality of
free topological groups \cite{Arh4}.
\item Tkachenko's result on the Souslin property of free
topological groups on compacta \cite{Tk2}.
\item Results on completeness and subgroups of free
topological groups by Uspenskij \cite{U9} and Sipacheva
\cite{Sip,Sip2}.
\item Results by Galindo and Hern\'andez on the reflexivity of
free abelian topological groups \cite{GH}.
\item Three of the results
reproduced below, namely the Tkachenko--Uspenskij Theorem 
\ref{tu} and Theorems \ref{dim} and \ref{lmp}.
\end{itemize}

The following result establishes a remarkable connection between
two of the objects so far introduced. It
was proved by successive efforts of
Tkachenko \cite{Tk}, who announced the result but supplied it
with a flawed proof, and Uspenskij \cite{U9}, who found both the flaw
and a correct proof some years later.

\begin{thm}[Tkachenko--Uspenskij, 1983/90]
For every pseudometric $\rho$ on the set $X^\dag=X\cup\{0\}$,
\[p_\rho\vert A(X) =\bar\rho.\]
As a corollary, for every topological space $X$ the free
abelian topological group $A(X)$ canonically embeds into the free
locally convex space $L(X)$ as a closed topological subgroup.
\qed
\label{tu}
\end{thm}

\begin{remark}
What is the noncommutative analogue of Tkachenko--Uspenskij Theorem?
Or, using the fashionable buzzword, how to {\it quantize} the above
result? \par
Firstly, we suggest that such a `quantization' must have to do with the
free SIN group $F_{\mathrm{SIN}}(X)$ rather than with the free topological
group $F(X)$ (which is really too complicated an object).\par
Also, notice that the Tkachenko--Uspenskij Theorem is equivalent
to the following statement: continuous homomorphisms from
the free abelian topological group $A(X)$ to the additive groups
of Banach spaces determine the topology of $A(X)$. It is not
difficult to verify that the additive topological group of every
Banach space embeds, as a topological subgroup, into the unitary
group of an abelian $C^\ast$-algebra, equipped with the induced
norm topology. A quantization of the above
statement will amount to allowing for all, and not just abelian, 
$C^\ast$-algebras. Finally notice 
that the unitary group of every $C^\ast$-algebra embeds, as a
topological subgroup,
into the unitary group of a Hilbert space equipped with the
uniform operator topology.
\par
Hence the resulting conjecture.
\end{remark}

\begin{conj}[`Non-commutative Tkachenko--Uspenskij Conjecture'] 
Continuous homomorphisms from the free balanced topological
group $F_{\mathrm{SIN}}(X)$ to the unitary groups $U(\H)_u$ of Hilbert
spaces with the uniform operator topology determine the topology
of $F_{\mathrm{SIN}}(X)$ for every (Tychonoff) topological space $X$. 
\label{nctuc}
\end{conj}

If the reader is unconvinced that the above is the `right' 
non-commutative version of Theorem \ref{tu}, notice that the
Tkachenko-Uspenskij theorem can be recast, modulo duality
theory for locally convex spaces, in the following equivalent form.
Let $I$ denote a fixed convex closed simply connected neighbourhood
of zero in the circle rotation group ${\mathbb{T}}=U(1)$.
Call a subset $A\subseteq G$ of an abelian topological
group $G$ a {\it polar set} if for some family $\mathcal X$ of continuous
characters of $G$ one has 
\[A=\bigcap_{\chi\in{\mathcal X}}\chi^{-1}(I).\]
Then Theorem \ref{tu} is equivalent to the statement that
the free abelian topological group $A(X)$ admits a neighbourhood
base consisting of polar sets. This was proved in \cite{Pdual}.

While proceeding from the abelian to non-abelian case, it is
natural to replace characters with finite-dimensional unitary
representations. Call a subset $A$ of a 
topological group $G$ a {\it polar set} if for some family
$\Pi=\cup_{n\in\N}\Pi_n$ of continuous finite-dimensional
unitary representations of $G$ there are simply connected,
convex (in the Riemannian sense) closed neighbourhoods of the
identity $I_\pi$ in the groups $U(n)$, where $\pi\colon G\to U(n)$,
such that 
\[A=\bigcap_{n\in\N}\bigcap_{\pi\in\Pi_n}\pi^{-1}(I_\pi).\]
It is an easy exercise to check
that every polar set is invariant.
Using the results on the so-called
residual finite-dimensionality of free $C^\ast$-algebras on
metric spaces \cite{Prfd}, one can prove that the 
non-commutative Tkachenko--Uspenskij conjecture
is equivalent to the existence of a neighbourhood base in the
topological group $F_{SIN}(X)$ consisting of polar sets in the above
sense.  
In such a form, the relationship
between Conjecture \ref{nctuc} and Theorem \ref{tu}
becomes obvious.

\subsubsection{}
Now we are fully armed to accomplish the proof of Theorem \ref{ab}. For
simplicity, we will only outline it in the case where the abelian
topological group
$G$ in question is not only countable, but metrizable as well;
the non-metrizable case only requires some extra
technical ingenuity but nothing really deep.
Denote by $\rho$ a translation-invariant metric generating the
topology of $G$. Since $G$ is clearly a topological factor-group
of the free group $A(G)$ equipped with the Graev metric $\bar\rho$,
it is enough to prove the theorem for the metric group 
$(A(G),\bar\rho)$ and then divide the monothetic group by
the kernel of the factor-homomorphism $A(G)\to G$: indeed,
monotheticity is preserved by proceeding to the images under
continuous homomorphisms with dense image. 
According to Tkachenko--Uspenskij Theorem, it suffices to prove
the statement for the separable normed space
$(L(G),p_\rho)$, which contains $(A(G),\bar\rho)$ as a topological
subgroup, or, equivalently, for the Banach space completion of
$(L(G),p_\rho)$, which we will denote for simplicity by $E$.
Enumerate by integers a countable everywhere dense subset
$\{x_n\colon n\in\N_+\}$ in $E$. Denote by $H$ the Banach space
direct sum of $E$ with the separable Hilbert space
$l_2(\N_+\times\N_+)$, and let $D$ be the subgroup of $H$
generated by all elements of the form
\[(mx_n, e_{m,n}),~ m,n\in\N_+.\]
(Here $e_{m,n}$ denote the standard basic vectors in the Hilbert
space $l_2(\N_+\times\N_+)$.)
Then one can verify that the Banach space $E$
is isomorphic to a topological subgroup of
the topological factor-group $H/D$, and the latter group
is clearly generated by the union of all its subgroups isomorphic
to the circle rotation group $\mathbb T$ (images of all
one-dimensional linear spaces passing through elements of $D$!)
and therefore $H/D$ is monothetic by the force of the following
result which we call Rolewicz Lemma. (Cf. \cite{Ro}, in which note
Rolewicz has actually proven the result below --- by an accurate,
recursive application of the Kronecker Lemma --- even if he never
stated the result in full generality and instead established it 
just for a concrete
example of a topological group he was constructing.)

\begin{thm}[Rolewicz Lemma]
A complete metric
abelian topological group $G$ topologically generated by
the union of countably infinitely many subgroups topologically
isomorphic to the circle group ${\Bbb T}\cong
U(1)$ is monothetic (that is, has one topological generator). \qed
\end{thm}

\subsection{Free groups: discrete {\it vs} topological}
In some respects, the known parallels between discrete and
topological groups go surprisingly far, especially for free groups.
Recall that two free bases in a free group always have the same
cardinality (called the {\it rank} of the free group). In other words,
if $X$ and $Y$ are two sets such that the free groups $F(X)$ and $F(Y)$
are isomorphic, then $\abs X=\abs Y$.
The same holds true for free abelian groups.
The following result can be seen as a perfect topological counterpart,
where the cardinality of a set is replaced with the dimension of a
topological space.

\begin{thm}[Pestov, 1982, \cite{P3}]
Let $X$ and $Y$ be two Tychonoff topological spaces with the property
that the free topological groups $F(X)$ and $F(Y)$ are 
isomorphic. {\rm (}Or: the free abelian topological groups,
$A(X)$ and $A(Y)$, are isomorphic.{\rm )} 
Then $X$ and $Y$ have the same Lebesgue covering
dimension: $\dim X=\dim Y$.
\label{dim}
\qed
\end{thm}

Parallels in this direction go further. It is known that the
free group $F_\infty$ on countably infinitely many generators 
embeds, as a subgroup, into the free group $F_2$ on two generators.
For topological groups, the following can be seen as a sensible
approximation to the same phenomenon.

\begin{thm}[Katz, Morris, Nickolas \cite{KMN}]
If $X$ is a countable CW-complex of finite dimension, then the
free topological group $F(X)$ is isomorphic with a topological
subgroup of $F({\mathbb I})$. \qed
\label{cm}
\end{thm}

What is the situation in the abelian case? The `discrete' suggestion
is that the rank of a subgroup of a free abelian group $A(X)$ cannot
exceed the rank of $A(X)$. For a while it remained unknown whether the
free abelian topological group on $A({\mathbb I}^2)$ embeds, as a
topological subgroup, into $A({\mathbb I})$. 
The answer turned out to be rather unexpected and requiring much 
subtler and advanced tools to obtain than for example the 
noncommutative theorem \ref{cm}. The following result demonstrates that
the analogy with the discrete case is, after all, not
comprehensive.

\begin{thm}[Leiderman--Morris--Pestov, \cite{LMP}]
If $X$ is a finite-dimensional metric compactum, then
$A(X)\hookrightarrow A(\I)$ as a topological subgroup.
\label{lmp}
\end{thm}

The proof is based on the following deep result, which in its time
had answered Hilbert's Problem 13. 

\begin{thm}[Kolmogorov Superposition Theorem, \cite{Kol, Ost}]
Every finite-dimen\-sional metric compactum $X$ possesses a
{\it basic system} of continuous
functions 
\[f_1,\dots,f_N\colon X\to {\mathbb I},\]
 meaning that every
continuous function $g\colon X\to {\mathbb I}$ can be represented as
the sum of compositions 
\[g=\sum_{i=1}^N h_i\circ f_i\]
of the basic functions with suitably chosen continuous functions
$h_i\colon {\mathbb I}\to {\mathbb I}$. \qed
\end{thm}

Let us show how to prove Theorem \ref{lmp}. 
First of all,
it is rather obvious that the space $C_p(\I)$ of continuous functions
with the pointwise topology admits a continuous linear operator onto
the space $C_p(\I\oplus\cdots\oplus \I)$ on the disjoint sum of finitely
copies of the interval. 
Using a basic system of functions on a finite-dimensional metrizable
compactum $X$, the space $C_p(\I\oplus\cdots\oplus \I)$ can be mapped
in a linear continuous surjective fashion onto the space $C_p(X)$.
An application of a theorem by Arhangel'ski\u\i\ \cite{Arh3}
enables one to conclude that the composition 
operator $C(\I)\to C(X)$ remains continuous with respect to the
compact-open (that is, uniform) topologies on both function spaces.
According to the Open Mapping Theorem, the continuous linear
operator $C(\I)\to C(X)$ is open. 
Duality theory for locally convex spaces leads one to invert the
direction of this operator and to obtain a topological embedding
$L(X)\hookrightarrow L(\I)$ of the free locally convex spaces.
Finally, one invokes
the Tkachenko-Uspenskij Theorem together with an observation that
the free abelian topological groups on $X$ and on $\I$ sit inside
the corresponding free locally convex spaces in the right
way to obtain a topological group embedding 
$A(X)\hookrightarrow A(\I)$. \qed

\begin{remark} Theorem \ref{lmp} suggests that problems about
free topological groups can be very difficult and substantial.
Imagine that someone has proved the
existence of a topological embedding $A(\I^2)\hookrightarrow A(\I)$
without prior knowledge of the Kolmogorov Superposition Theorem;
such a person would be forced to essentially rediscover the solution
to Hilbert's Problem 13 on his/her own. 
Perhaps something of the kind indeed happens in general
topology from time to time,
and if anything, this shows the need for us the general topological
algebraists to consciously look out for links
with other areas of mathematics.
\end{remark}

\begin{remark} And not just mathematics.
There is, it seems, an interesting perspective of linking
theory of free topological groups to computing.

Let $X=(X,\rho)$ be a metric space, where the value of the metric
$\rho(x,y)$ is interpreted as the {\it cost} 
of transporting a unit mass from point $x$ to point $y$.
Suppose a unit mass is distributed between and
stored at points $x_1,x_2,\cdots, x_n$,
with amounts $\lambda_1,\lambda_2,\cdots,\lambda_n$ at each of them,
and we want to transport the mass and store it at points
$y_1,y_2,\cdots,y_m$, with amount $\mu_1,\mu_2,\cdots,\mu_n$ at
each of the corresponding points. The cost of performing such 
a transportation is known as the {\it Kantorovich distance}
between $\sum_{i=1}^n\lambda_i x_i$ and
$\sum_{j=1}^n\mu_jy_j$. The Kantorovich distance plays a singularly
important role in a wide range of applied mathematical sciences,
from probability theory through information theory to computing and
data storage and analysis (for a recent comprehensive treatment, see
the two-volume set \cite{RR}
exclusively devoted to the Kantorovich distance).
At the same time, it is not difficult to see that the Kantorovich
distance is exactly the metric generated by the maximal norm $p_\rho$ on
the free locally convex space. As such, it can be 
approximated with any given degree of accuracy by the
Graev metric on $A(X)$. 

The Kantorovich distance can be computed through linear programming
in quadratic time in the input size, $n$.
An outstanding problem about the Kantorovich distance is the following:
does there exist an algorithm for computing the value of the
distance in the {\it linear} time $O(n)$ in the size of the input, $n$? 
Currently such algorithms are only known for $X=\R$ with the usual
distance and also $X=\s^1$, the circle. The
accumulated combinatorial techniques for dealing with the
Graev metric make free topological group theorists well-poised to tackle
this problem  for more general metric spaces than the
real line, especially in the natural 
case where $X$ is a finite-dimensional
compactum.
For example, can the embedding 
$A(X)\hookrightarrow A({\mathbb I})$ constructed in the proof of
Theorem \ref{lmp} be used to
achieve an algorithmic speed-up?

In the non-abelian case, the Graev metric on the free group $F(X)$
can be shown to coincide with another distance of great importance
in computing: the so-called {\it string edit distance}
between finite words in an alphabet $W$. 
The value of the string edit distance between two strings
(words), $\sigma$ and $\tau$, is defined as the 
minimal number of insertions, deletions, and replacements necessary
to get $\sigma$ from $\tau$.  (Example:
$d({\tt metric, track})=4$.) Clearly, the string edit distance is
recovered from the Graev extension of the discrete metric from $W$
over $F(W)$. This type of distance is of importance in particular
in molecular biology \cite{AG}, 
and it is not impossible that the fluency
of free topological group theorists in manipulating Graev metrics
can be used to improve on the existing algorithms for computing
the string edit distance and its various modifications.
\end{remark}

\subsection{The epimorphism problem} 
The topic with which I would like to conclude the article provides
a nice example where all our main lines of development
--- actions on compacta, `massive' groups, free topological groups,
and embeddings --- 
converge and work together in harmony.

A morphism $f\colon H\to G$ between two objects in some category is called
an {\it epimorphism} if for all objects $F$ and two arbitrary morphisms 
$g,h\colon G\to F$, the condition
$g\circ f=h\circ f$, cf. the diagram
\[H  \stackrel{f}{ \to} G \stackrel{g}
{\rightrightarrows}_{_{\!\!\!\!\!\! h}} F,\]
implies that $g=h$.

Let us first consider the category of discrete groups and group
homomorphisms. Certainly, every homomorphism onto is an epimorphism.
It turns out that the converse is also true. The following is apparently
a part of the group theory folklore.

\begin{thm}
A group homomorphism $f\colon H\to G$ is an epimorphism if and only
if $f$ is onto.
\label{discr}
\end{thm}

It is both illuminating and essential for what follows to go through
the proof. Let $f$, $H$ and $G$ be as in the statement
of the theorem. Assume that $f$ is not onto, and let us show that $f$
is not an epimorphism. To produce a group $F$ and two distinct morphisms
$g,h\colon G \to F$ with the property $g\circ f = h\circ f$, we proceed
as follows. The image $f(H)$ is a proper subgroup of $G$.
Denote by $X=G/f(H)$ the left factor-set, that is, the
collection of all left cosets $x\cdot f(H)$, where $x\in G$. 
Then $\abs X>1$. 
The natural action of $G$ on $X$ by left translations extends to an
action of $G$ on the free group $F(X)$ by group isomorphisms,
that is, every element $x\in G$ determines a group isomorphism
$y\mapsto x\cdot y$ of $F(X)$, and the resulting mapping
from $G$ to the group of automorphisms of $F(X)$ is a homomorphism
of groups.
Now one can form the {\it semidirect product},
$G\ltimes F(X)$, corresponding to such an action. This is a group,
which is, set-theoretically, the Cartesian product $G\times F(X)$,
equipped with the group operation as follows:
\[(x,y)(x',y') := (xx', y(x\cdot y')).\]
The identity is the element $(e_G, e_{F(X)})$, while the inverse of
an element $(x,y)$ is simply $(x^{-1},x^{-1}\cdot (y^{-1}))$. 
Notice that $G$ forms
a subgroup (and at the same time a factor-group) of the semidirect
product, under the embedding $x\mapsto (x,e)$. 
Set
$F=G\ltimes F(X)$, and define 
two homomorphisms from $G$ to $F$ 
by letting $g$ be the above described embedding
$G\hookrightarrow G\ltimes F(X)$, while $h$ sends an element $x\in G$
to its conjugate in $F$ by the element $(e_G,f(H))$, where the subgroup
$f(H)$ is viewed as a coset and element of the free basis $X\subset F(X)$.
It is easy to see that $g\vert_H=h\vert_H$, while in general the
two homomorphisms $g$ and $h$ differ
(in fact, $g(x)\neq h(x)$ whenever $x\in G\setminus f(H)$). The proof is
finished.

What happens in the category of (Hausdorff! -- here it is really
essential) topological groups and continuous
homomorphisms? The most immediate observation is that if
$f\colon H\to G$ is a continuous homomorphism with everywhere dense
range ($\overline{g(H)}=G$), then $f$ is an epimorphism.
It was asked by Karl H. Hofmann in late 1960's
whether the converse was true. In other words, must
every epimorphism between two topological groups have a dense range?
As we have just seen, the answer is
`yes' if $G$ is discrete. It is an obvious exercise to show that
the answer is positive if $G$ is abelian. By the end of this section,
the reader will be able to prove that the answer is positive if
$G$ is locally compact.\footnote{Using the following
result by de Vries \cite{dV}: 
every continuous action of a locally compact
group $G$ on a topological space $X$ is {\it linearizable,} that is,
there exists a continuous action of $G$ on a locally convex space
$E$ by isomorphisms and an embedding $X\hookrightarrow E$ as
a topological $G$-subspace.}
The general case, however, was
only settled (in the negative) a few years ago by means of the following
astonishing result. (Cf. \cite{U7,U8} for further refinements.)

\begin{thm}[Uspenskij, 1993, \cite{U6}]
Let $x$ be any point of the circle $\s^1$. Then the embedding
$\mathrm{St}_x\hookrightarrow\Homeo({\mathbb S}^1)$ 
is an epimorphism, where $\mathrm{St}_x=\{g\vert g\cdot x=x\}$ 
is the isotropy subgroup.
\label{isot}
\end{thm}

To understand how such a seemingly improbable thing 
can happen, it is worth re-examining the proof 
of Theorem \ref{discr} and finding out to what extent one can have
it `topologized.' \par
It poses no problem at all, to make the semidirect
product of two topological groups into a topological group: the
usual product topology will do the trick, provided the action of
the first group on the second is continuous, as a map (in our case)
\begin{equation}
G\times F(X)\to F(X).
\label{action}
\end{equation}
(The standard reference for semidirect 
products of topological groups is \cite{HR}, 2.6.20.)
The factor-space $X=G/\overline{f(H)}$ is of course a topological
space, so that it is natural to equip $F(X)$
with the topology of the free topological group
on $X$. Here we approach
the bottleneck of the argument: the action (\ref{action}) need not be
continuous! What is worse, the situation cannot be, in general, remedied
by considering group topologies on $F(X)$ other than the free
topology. A convenient framework enabling one
to deal with situations of this kind had been developed by
Megrelishvili independently from Uspenskij's Theorem \ref{isot} 
(though published later, \cite{Meg4, Meg5}),
and we will now describe it briefly. 

Let $G$ and $F$ be two topological groups.
Say that $F$ is a $G${\it -group} if  
$G$ continuously acts on $F$ by automorphisms. More precisely,
there is given a continuous action $G\times F\to F$, and every
motion $F\ni y\mapsto x\cdot y\in F$ is a group automorphism. 
(A particular case of this situation arises whenever a topological
group $G$ is represented in a normed space $F$.) It is clear how to
define morphisms between two $G$-groups: they are continuous
homomorphisms commuting with the action of $G$. 
Now, given a topological group $G$ acting continuously on a
topological space $X$, one can define the {\it free topological
$G$-group} on $X$ in a standard fashion. Namely,
$F_G(X)$ is a topological $G$-group, and there is 
an (essentially unique) morphism of $G$-spaces
$\iota\colon X\to F_G(X)$ with the property that each morphism of
$G$-spaces $f$ from $X$ to an arbitrary $G$-group $A$
admits a unique factorization $f=\bar f\circ\iota$, where
$\bar f\colon F_G(X)\to A$ is a morphism between $G$-groups.
For example, if $G=\{e\}$ is a  trivial group, 
then $F_G(X)$ turns into the free topological group on $X$.
The free topological $G$-group $F_G(X)$ is algebraically generated by
the set $\iota(X)$. For a given $G$-space $X$ the free
topological $G$-group $F_G(X)$ is 
the Hausdorff replica of the
free group $F(X)$ equipped with the finest group topology such that
its restriction to $X$ is coarser than the original topology on $X$
and the action $G\times F(X)\to F(X)$ is continuous.
One says that the free topological $G$-group is {\it trivial} if
it is isomorphic to $\Z$ equipped with the discrete topology and 
the trivial action of $G$, in which case $\iota$ is constant
and takes all of $X$ to $1\in\Z$. 

Getting back to our discussion of the proof of Theorem \ref{discr},
we can see that the right substitute for the free group
$F(X)$ is the free topological $G$-group $F_G(X)$ on the $G$-space
$X$. Therefore, everything boils down to the question of
nontriviality of $F_G(X)$. We have practically 
established the following.

\begin{thm}[Pestov, \cite{P18}]
A continuous homomorphism $f$ between topological groups $H$ and $G$
is an epimorphism if and only if the free topological $G$-group
on the $G$-space $X=G/\overline{f(H)}$ is trivial. \qed
\end{thm}

The {\it Effros microtransitivity theorem} \cite{Eff}
says that if a
Polish (= completely metrizable second-countable) topological group $G$
acts transitively on a Polish space $X$, then $X$ is isomorpic,
as a $G$-space, with the left factor-space of $G$ by the isotropy
subgroup $\mathrm{St}_x$ of an arbitrary element $x\in X$.
Now we obtain the following convenient corollary.

\begin{corol}[\cite{P18}]
 Let $X$ be a
transitive $G$-space such that both the acting group $G$ and the  
space $X$ are Polish. The following conditions are
equivalent:
\smallskip
\item{(i)} the free topological $G$-group $F_G(X)$ is trivial;
\item{(ii)} the canonical embedding of the isotropy subgroup
$\operatorname{St}_x$ of any point $x\in X$ into $G$ is
an epimorphism of Hausdorff topological groups.
\label{staa}
\qed
\end{corol}

The following result was circulated by Megrelishvili in a preprint form
in early 1990's.

\begin{thm}[Megrelishvili, \cite{Meg4, Meg5}]
The free topological $\Homeo({\mathbb I})$-group \\
$F_{\Homeo(\I)}(\I)$ is trivial. \qed
\end{thm}

As an obvious corollary, the free topological 
$\Homeo({\s^1})$-group $F_{\Homeo(\s^1)}(\s^1)$ is trivial as well,
and by invoking Corollary \ref{staa}, we obtain a proof of
Uspenskij's Theorem \ref{isot}.

\bibliographystyle{amsplain}

\end{document}